\documentclass[11pt]{amsart}
\usepackage{amsmath,amsthm,amssymb}
\usepackage[all]{xy}

\begin{document}

\newtheorem{theorem}{Theorem}[section]
\newtheorem{lemma}[theorem]{Lemma}
\newtheorem{proposition}[theorem]{Proposition}
\newtheorem{corollary}[theorem]{Corollary}
\newtheorem*{corollarya1}{Corollary A.1}
\newtheorem*{corollarya2}{Corollary A.2}
\newtheorem*{propositiona1}{Proposition A.1}
\newtheorem*{theoremb1}{Theorem B.1}
\newtheorem{bigtheorem}{Theorem}

\theoremstyle{definition}
\newtheorem{definition}[theorem]{Definition}
\newtheorem{example}[theorem]{Example}
\newtheorem{formula}[theorem]{Formula}
\newtheorem{nothing}[theorem]{}

\theoremstyle{remark}
\newtheorem{remark}[theorem]{Remark}

\renewcommand{\arraystretch}{1.2}

\newcommand{\dual}{^{\vee}}
\newcommand{\contr}{{\mspace{1mu}\lrcorner\mspace{1.5mu}}}
\newcommand{\de}{\partial}
\newcommand{\debar}{{\overline{\partial}}}
\newcommand{\pibar}{\overline{\pi}}

\newcommand{\desude}[2]{{\dfrac{\de #1}{\de #2}}}

\newcommand{\mapor}[1]{{\stackrel{#1}{\longrightarrow}}}
\newcommand{\ormap}[1]{{\stackrel{#1}{\longleftarrow}}}

\newcommand{\mapver}[1]{\Big\downarrow\vcenter{\rlap{$\scriptstyle#1$}}}

\newcommand{\binfty}{\boldsymbol{\infty}}
\newcommand{\bi}{\boldsymbol{i}}
\newcommand{\bl}{\boldsymbol{l}}

\renewcommand{\bar}{\overline}
\renewcommand{\Hat}[1]{\widehat{#1}}

\newcommand{\sA}{\mathcal{A}}
\newcommand{\Oh}{\mathcal{O}}
\newcommand{\sH}{\mathcal{H}}
\newcommand{\sL}{\mathcal{L}}
\newcommand{\sM}{\mathcal{M}}
\newcommand{\sB}{\mathcal{B}}
\newcommand{\sY}{\mathcal{Y}}

\newcommand{\Q}{\mathbb{Q}}
\newcommand{\K}{\mathbb{K}}
\newcommand{\Proj}{\mathbb{P}}

\newcommand{\DER}{{{\mathcal D}er}}

\newcommand{\ad}{\operatorname{ad}}
\newcommand{\MC}{\operatorname{MC}}
\newcommand{\Def}{\operatorname{Def}}
\newcommand{\Hom}{\operatorname{Hom}}
\newcommand{\End}{\operatorname{End}}
\newcommand{\Image}{\operatorname{Im}}
\newcommand{\image}{\operatorname{Im}}
\newcommand{\Der}{\operatorname{Der}}
\newcommand{\Mor}{\operatorname{Mor}}
\newcommand{\Cone}{\operatorname{Cone}}
\newcommand{\Aut}{\operatorname{Aut}}
\newcommand{\coker}{\operatorname{Coker}}

\newcommand{\Grass}{\operatorname{Grass}}
\newcommand{\Spec}{\operatorname{Spec}}


\title[$L_{\infty}$-algebras, Cartan homotopies and period maps]
{$\mathbf{L}_{\boldsymbol\infty}$-algebras, Cartan homotopies and period maps}
\date{May 9, 2006}

\author{Domenico Fiorenza}
\address{\newline Dipartimento di Matematica ``Guido
Castelnuovo'',\hfill\newline
Universit\`a di Roma ``La Sapienza'',\hfill\newline
P.le Aldo Moro 5,
I-00185 Roma Italy.\hfill\newline}
\email{fiorenza@mat.uniroma1.it}
\urladdr{www.mat.uniroma1.it/\~{}fiorenza/}
\author{Marco Manetti}
\email{manetti@mat.uniroma1.it}
\urladdr{www.mat.uniroma1.it/people/manetti/}

\begin{abstract}
We prove that, for every compact
K\"{a}hler manifold, the period map of its Kuranishi family is induced by a
natural $L_{\infty}$-morphism.
This implies, by standard facts about $L_{\infty}$-algebras, that
the period map is a ``morphism of deformation theories'' and
then commutes with all
deformation theoretic constructions (e.g. obstructions).
\end{abstract}

\subjclass{14D07, 17B70, 13D10}
\keywords{Differential graded Lie algebras, symmetric coalgebras,
$L_{\infty}$-algebras, functors of Artin rings, K\"{a}hler manifolds, period map}

\maketitle

\tableofcontents

\vfill\eject
\section*{Introduction}

The homotopy Lie algebra approach to deformation theory over a
field $\K$ of characteristic 0
is based on two creeds.\\ First creed: let $\mathcal{M}$ be a moduli space
(for some
classification problem) and
$x$ a point of $\mathcal{M}$; then the geometry of  the formal
neighbourhood of $x$ in $\mathcal{M}$ is encoded into an (homotopy
class of) $L_\infty$-algebra(s).
More precisely, there exists an
$L_\infty$-algebra $\mathfrak{g}$, defined up to quasiisomorphism,
such that, for any local Artinian
$\K$-algebras $B$ with maximal ideal $\mathfrak{m}_B$ and
residue field $\K$ one has
\[\{\phi\colon \Spec(B)\to \mathcal{M}\mid
\phi(\Spec(\mathbb{K}))=x\}=\Def_\mathfrak{g}(B),\]
where
\[\Def_\mathfrak{g}(B)=\frac{\{\text{Solutions of the Maurer-Cartan equation in }
\mathfrak{g}\otimes_{\K}\mathfrak{m}_B\}}{\text{homotopy equivalence}}.\]
Second creed: every ``natural''
morphism between formal pointed moduli spaces is
induced by an $L_{\infty}$-morphism between the associated
$L_\infty$-algebras.\\
In this paper we shall explicate the first creed for
Grassmannians using a  general construction that
also applies to  other moduli spaces
(dg-Grassmannians, Quot and Hilbert schemes,  Brill-Noether loci etc.)
and the second creed for the universal period map,
intended as a natural morphism of moduli spaces:
from deformations of a compact
K\"ahler manifold to deformations of the Hodge filtration of
its De Rham cohomology.\\

More concretely, let's denote by
\[\mathcal{X}\to (S,0),\qquad  \mathcal{X}=\bigcup_{t\in S}X_t,\]
  the  Kuranishi
family of a compact K\"ahler manifold $X=X_0$;
let $p\ge 0$ be a fixed integer and consider
the period map \cite[10.1.2]{Voisin}
\[ \mathcal{P}^p\colon (S,0)\to \Grass(H^*(X,\mathbb{C}))=\prod_i 
\Grass(H^i(X,\mathbb{C})),\]
\[ t\mapsto F^pH^*(X_t,\mathbb{C})=\prod_i F^pH^i(X_t,\mathbb{C}).\]

Griffiths proved \cite{Griffiths} that $\mathcal{P}^p$ is a holomorphic map
and its differential $d\mathcal{P}^p$ is the same of
the contraction map
\[ \bi\colon H^1(X,T_X)\to
\bigoplus_i \Hom\left(F^pH^i(X,\mathbb{C}),
\frac{H^{i}(X,\mathbb{C})}{F^pH^{i}(X,\mathbb{C})}\right),\qquad
\bi_{\xi}(\omega)=\xi\contr\omega.\]
It is also known \cite{clemens,CCK,ranUVHS} that obstructions to
deformations of $X$ are contained in the kernel of
\[ \bi\colon H^2(X,T_X)\to
\bigoplus_i \Hom\left(F^pH^i(X,\mathbb{C}),
\frac{H^{i+1}(X,\mathbb{C})}{F^pH^{i+1}(X,\mathbb{C})}\right),\qquad
\bi_{\xi}(\omega)=\xi\contr\omega\]
(this fact is known as Kodaira's Principle \emph{ambient cohomology 
annihilates obstruction}).
However  the proofs of \cite{clemens,CCK,ranUVHS}
are not completely satisfying because the period map plays only a
marginal role in them, while the most natural  way
of proving Kodaira's Principle would be to show that the
period map is a ``morphism of deformation theories'', i.e.
that $\mathcal{P}^p$ commutes
with every deformation theoretic construction:
for instance obstruction theories.\\
It is well known that the deformations of a compact complex
manifold $X$ are governed by the $L_{\infty}$-algebra underlying the
Kodaira-Spencer algebra $K_X:=\oplus_{i} A^{0,i}_X(T_X)$.\\

The main results of this paper are:
\begin{enumerate}

\item We explicitly describe  an $L_{\infty}$-algebra $C^p$ such that
\[\{\phi\colon \Spec(B)\to \Grass(H^*(X,\mathbb{C}))\mid
\phi(\Spec(\mathbb{C}))=F^pH^*(X,\mathbb{C})\}=\Def_{C^p}(B).\]

\item We explicitly describe a linear $L_{\infty}$-morphism $K_X\to C^p$ inducing
$\mathcal{P}^p$ at the level of deformation functors.
\end{enumerate}

These results give us an algebraic description of the period map and imply,
by  general theory of $L_{\infty}$-algebras,
that $\mathcal{P}^p$ is a morphism of deformation theories.\\

The paper is divided in three parts: roughly speaking, in the first we 
make a functorial construction of an 
$L_{\infty}$ structure on the mapping cone of a morphism of differential
graded Lie algebras and we introduce the notion of Cartan homotopy.
A slightly expanded version of this part is available as \cite{cone}.\\
In the second part we describe an $L_{\infty}$-algebra governing deformations 
of subcomplexes and then giving a local description of Grassmannians.\\
Finally in the third part we exhibit an $L_{\infty}$-morphism inducing the universal 
period map of a compact K\"{a}hler manifolds.

\medskip
\textbf{Acknowledgment.} Our thanks to Jim Stasheff for precious
comments on the first version of \cite{cone}
and to Barbara Fantechi for useful discussions about dg-Grassmannians.\\

\textbf{Keywords and general notation.}

We assume that the reader is familiar with the notion and main
properties
of differential graded Lie algebras and
$L_{\infty}$-algebras (we refer to
\cite{fuka,grassi1,K,LadaMarkl,LadaStas,defomanifolds}
as introduction of such structures); however
the basic definitions are recalled in this paper in order to fix
notation and terminology.\\
For the whole paper, $\mathbb{K}$
is a field of characteristic 0; every vector space
is intended  over
$\mathbb{K}$.\\
$\mathbf{Art}$ is the category of local Artinian
$\K$-algebras with residue field $\K$.
For $A\in\mathbf{Art}$ we denote by
$\mathfrak{m}_A$ the maximal ideal of $A$.

\bigskip
\part{$L_{\infty}$ structures on mapping cones}

Let $\chi\colon L\to M$ be a morphism of differential graded Lie
algebras over a field $\K$ of
characteristic 0.
In the paper \cite{semireg} one of the authors has introduced, having
in mind
the example of embedded
deformations, the notion of Maurer-Cartan equation and gauge action
for the
triple
$(L,M,\chi)$; these notions reduce to the standard
Maurer-Cartan equation and gauge action of $L$ when $M=0$.
More precisely there are defined two functors of Artin rings
$\MC_{\chi},\Def_{\chi}\colon \mathbf{Art}\to \mathbf{Set}$,
where $\mathbf{Art}$ is the category of local Artinian $\K$-algebras
with residue field $\K$, in the
following way:
\[ \MC_{\chi}(A)=
\left\{(x,e^a)\in (L^1\otimes\mathfrak{m}_A)\times \exp(M^0\otimes
\mathfrak{m}_A)\mid
dx+\frac{1}{2}[x,x]=0,\;e^a\ast\chi(x)=0\right\},\]
where for $a\in M^0\otimes
\mathfrak{m}_A$ we denote by $\ad_a\colon M\otimes\mathfrak{m}_A\to
M\otimes\mathfrak{m}_A$
the operator $\ad_a(y)=[a,y]$
and
\[ e^a\ast
y=y+\sum_{n=0}^{+\infty}\frac{\ad^n_a}{(n+1)!}([a,y]-da)
=y+\frac{e^{\ad_a}-1}{\ad_a}([a,y]-da),\qquad y\in
M^1\otimes\mathfrak{m}_A.\]%
Then one defines
\[ \Def_\chi(A)=\frac{\MC_{\chi}(A)}{\text{gauge equivalence}},\]
where two solutions of the Maurer-Cartan equation are gauge
equivalent
if they belong to the same orbit of the gauge action
\[(\exp(L^0\otimes\mathfrak{m}_A)\times
\exp(dM^{-1}\otimes\mathfrak{m}_A))
\times \MC_{\chi}(A)\mapor{\ast}\MC_{\chi}(A)\]%
given by the formula
\[ (e^l, e^{dm})\ast (x,e^a)=(e^l\ast x,
e^{dm}e^ae^{-\chi(l)})=(e^l\ast x,
e^{dm\bullet a\bullet (-\chi(l))}).\]%
The $\bullet$ in the rightmost term in the above formula is
the Baker-Campbell-Hausdorff multiplication; namely
$e^xe^y=e^{x\bullet y}$.\\
Several examples in \cite{semireg} illustrate the utility of this
construction
in deformation theory.\\
In the same paper it is also proved that $\Def_\chi$ is the
truncation of an
extended  deformation functor \cite[Def. 2.1]{EDF} $F$ such that
$T^iF=H^i(C_{\chi})$,
where $C_{\chi}$ is the differential graded vector space
$C_\chi=(\mathop{\oplus}_iC_\chi^i,\delta)$,
\[C_\chi^i=
L^i\oplus M^{i-1},\qquad \delta(l,m)=(dl,\chi(l)-dm),\qquad l\in
L,m\in M.\]
By a general result \cite[Thm. 7.1]{EDF} there exists an $L_{\infty}$
structure on $C_{\chi}$,
defined up to homotopy,
whose associated deformation functor is isomorphic
to $\Def_{\chi}$.\\
The main result of this part is to describe explicitly a canonical
(and hence functorial)
$L_{\infty}$ structure on $C_{\chi}$ with the above property; this is
done in an elementary way, without using
the theory of extended deformation functors, so that the result of this part
can be used to reprove the main results of \cite{semireg}
without using \cite{EDF}.

\bigskip

\section{Conventions on graded vector
spaces}\label{sec.graded.conv.}
In this paper we will work with $\mathbb{Z}$-graded vector spaces;
we write a graded vector space as $V=\oplus_{n\in {\mathbb Z}}V^n$,
and call $V^n$ the degree $n$ component of $V$; an element $v$ of
$V^n$ is called a degree $n$ homogeneous element of $V$.

We adopt the convention according to which
degrees are `shifted on the left'. By this we mean that, for every integer $n$,  
$V[n]\simeq
\mathbb{K}[n]\otimes V$ where $\mathbb{K}[n]$ denotes the graded vector
space consisting in the field $\mathbb{K}$ concentrated in
degree $-n$. Note that, with this convention the canonical
isomorphism $V\otimes \mathbb{K}[n]\simeq V[n]$ is
$v\otimes 1_{[n]}\mapsto (-1)^{n\deg(v)}v_{[n]}$ and we have the
following isomorphism, usually called \emph{decalage} 
\begin{align*}
V_1[1]\otimes\cdots \otimes V_n[1]&\xrightarrow{\sim}
(V_1\otimes\cdots \otimes V_n)[n]\\
{v_1}_{[1]}\otimes \cdots \otimes{v_n}_{[1]}&\mapsto
(-1)^{\sum_{i=1}^n(n-i)\cdot\deg{v_i}}(v_1\otimes \cdots \otimes
v_n)_{[n]}.
\end{align*}

Denote by $\bigotimes^n V$, $\bigodot^n V$ and $\bigwedge^n V$
the $n$-th tensor, symmetric and exterior powers of $V$ respectively.
As with ordinary vector spaces, one can identify $\bigodot^nV$ and
$\bigwedge^n V$ with suitable subspaces of $\bigotimes^nV$, called the
subspace of symmetric and antisymmetric tensors respectively.
The decalage induces a canonical isomorphism
\[
\bigodot^n(V[1])\xrightarrow{\sim}\left(\bigwedge^n V\right)[n].
\]

\begin{remark}\label{rem.suspension} 
Using the natural isomorphisms
\[ \Hom^i(V,W[l])\simeq \Hom^{i+l}(V,W)\]
and the decalage isomorphism, we obtain natural identifications
\[ \operatorname{dec}\colon\Hom^i\left(\bigotimes^k V,W\right)\xrightarrow{\sim}  
\Hom^{i+k-l}\left(\bigotimes^k(V[1]),W[l]\right),\]
where
\[ \operatorname{dec}(f)(v_{1[1]}\otimes \cdots\otimes v_{k[1]})=
(-1)^{ki+\sum_{j=1}^k(k-j)\cdot\deg(v_j)}f(v_1\otimes \cdots\otimes v_{k})_{[l]}.\]
By the above considerations
\[ \operatorname{dec}\colon\Hom^i\left(\bigwedge^k V,V\right)\xrightarrow{\sim}  
\Hom^{i+k-1}\left(\bigodot^k(V[1]),V[1]\right).\]

\end{remark}

\bigskip
\section{Differential graded Lie algebras and $L_\infty$-algebras}
\label{sec.linfty.intro}

A differential graded Lie algebra (DGLA for short) is a Lie algebra
in the category of graded vector spaces, endowed with a compatible degree 1 differential.
More explicitly, it is the
data $(V,d,[\,,\,])$, where $V$ is a graded vector space, the Lie
bracket
\[
[\,,\,]\colon V\wedge V\to V
\]
satisfies the graded Jacobi identity:
\[
[v_1,[v_2,v_3]]=[[v_1,v_2],v_3]+(-1)^{\deg(v_1)\deg(v_2)}[v_2,
[v_1,v_3]],
\]
and  $d\colon
V\to V$ is a degree 1 differential, i.e.,
\[ d^2=0,\qquad
d[v_1,v_2]=[dv_1,v_2]+(-1)^{\deg(v_1)}[v_1,dv_2].
\]
Via the
decalage isomorphisms one can look at the Lie bracket of a DGLA $V$
as  a morphism
\[
q_2^{}\in\Hom^1(V[1]\odot V[1],V[1]),\qquad q_2(v_{[1]}\odot
w_{[1]})=(-1)^{\deg(v)}[v,w]_{[1]},
\]
Similarly, the suspended differential $q_1=d_{[1]}={\rm
id}_{\K[1]}\otimes d$ is a morphism of degree 1
\[
q_1^{}\colon V[1]\to V[1],\qquad q_1(v_{[1]})=-(dv)_{[1]}.
\]
Up to the canonical bijective linear map $V\to V[1]$, $v\mapsto v_{[1]}$,
the suspended differential $q_1$ and the bilinear operation
$q_2$ are written simply as
\[
q_1(v)=-dv; \qquad
q_2(v\odot w)=(-1)^{\deg_V(v)}[v,w].
\]
   Define morphisms $q_k^{}\in \Hom^1(
\odot^k(V[1]), V[1])$ by setting $q_k^{}\equiv 0$, for
$k\geq 3$. The map
\[
\sum_{n\geq 1} q_n\colon \bigoplus_{n\geq 1}\bigodot^n V[1]\to
V[1]
\]
extends to a coderivation of degree 1
\[ Q\colon \bigoplus_{n\ge 1}\bigodot^n V[1]\to \left(\bigoplus_{n\ge
1}\bigodot^n
V[1]\right)\]
on the reduced symmetric coalgebra cogenerated by $V[1]$, by the
formula
\begin{nothing}\label{not.codifferential}
\[Q(v_1\odot\cdots\odot v_n)=\sum_{k=1}^n\sum_{\sigma\in S(k,n-k)}
\varepsilon(\sigma)q_k(v_{\sigma(1)}\odot\cdots\odot v_{\sigma(k)})
\odot v_{\sigma(k+1)}\odot\cdots\odot v_{\sigma(n)},\]
\end{nothing}
where
$S(k,n-k)$ is the set of unshuffles and
$\varepsilon(\sigma)=\pm 1$ is the  \emph{Koszul sign},
determined by the relation in
$\bigodot^n V[1]$
\[ v_{\sigma(1)}\odot\cdots\odot v_{\sigma(n)}=\varepsilon(\sigma)
v_{1}\odot\cdots\odot v_{n}.\]
The axioms of differential graded Lie algebra are then equivalent
to $Q$ being a codifferential, i.e., $Q^2=0$. This description of
differential graded Lie algebras in terms of the codifferential
$Q$ is called the Quillen construction \cite{Qui}. By
dropping the requirement that
$q_k^{}\equiv 0$ for
$k\geq 3$ one obtains the notion of $L_\infty$-algebra (or strong
homotopy Lie algebra), see e.g. \cite{LadaMarkl,LadaStas,K};
namely, an
$L_{\infty}$ structure on a graded vector space
$V$ is a sequence
of linear maps of degree 1
\[ q_k\colon \bigodot^k V[1]\to
V[1],\qquad k\ge 1,\]
such that the induced coderivation $Q$
on the reduced symmetric coalgebra cogenerated by $V[1]$, given  by
the Formula~\ref{not.codifferential} is a codifferential, i.e.
i.e., $Q^2=0$. This condition in particular implies $q_1^2=0$, i.e.,
an $L_\infty$-algebra is in particular a differential complex. Note
that, by the above
discussion, every DGLA can be naturally seen as an
$L_\infty$-algebra; namely, a DGLA is an $L_\infty$-algebra with
vanishing higher multiplications $q_k^{}$, $k\geq 3$. Note that, via the
decalage isomorphisms of Remark~\ref{rem.suspension}, 
every component $q_k^{}$ of an
$L_\infty$ structure on  $V$ can be seen as morphism
\[
\mu_k^{}\in\Hom^{2-k}(\bigwedge^k V, V).
\]

\par A morphism
$f_\infty^{}$ between two $L_\infty$-algebras
$(V,q_1^{},q_2{},q_3^{},\dots)$ and  $(W,p_1^{},p_2{},p_3^{},\dots)$
is a sequence of linear
maps of degree 0
\[ f_n\colon \bigodot^n V[1]\to
W[1],\qquad n\ge 1,\]
such that the morphism of coalgebras
\[ F\colon \bigoplus_{n\ge 1}\bigodot^n V[1]\to \bigoplus_{n\ge
1}\bigodot^n
W[1]\]
induced by $\sum_n f_n\colon\bigoplus_{n\ge 1}\bigodot^n V[1]\to
W[1]$
commutes with the codifferentials induced by the two $L_{\infty}$
structures on $V$ and
$W$ \cite{fuka,K,LadaMarkl,LadaStas,defomanifolds}.
An $L_{\infty}$-morphism $f^{}_\infty$ is called
\emph{linear}  (sometimes \emph{strict}) if $f_n=0$ for every
$n\ge 2$.  We note that a linear map $f_1\colon V[1]\to W[1]$
is a linear $L_{\infty}$-morphism if and only if
\[ p_n(f_1(v_1)\odot\cdots\odot f_1(v_n))=f_1(q_n(v_1\odot\cdots\odot
v_n)),\qquad
\forall\; n\ge 1,\; v_1,\ldots,v_n\in V[1].\]
For instance, 
morphisms between DGLAs are
linear morphisms between the corresponding $L_\infty$-algebras.
\par
If $f_\infty$ is an $L_\infty$-morphism between
$(V,q_1^{},q_2{},q_3^{},\dots)$ and  $(W,p_1^{},p_2{},p_3^{},\dots)$,
then its linear part
\[
f_1\colon V[1]\to W[1]
\]
satisfies the equation $f_1\circ q_1=p_1\circ f_1$, i.e., $f_1$ is a
map of differential complexes $(V[1],q_1)\to (W[1],p_1)$. An
$L_\infty$-morphism $f^{}_\infty$ is called a quasiisomorphism of
$L_\infty$-algebras if its linear part $f_1$ is a quasiisomorphism of
differential
complexes. A major result in the theory of $L_\infty$-algebra is the
following \emph{homotopical transfer of structure} theorem (see
\cite{fuka,KonSoi} for a proof).
\begin{theorem}
Let $(V,q_1^{},q_2{},q_3^{},\dots)$ be an $L_\infty$-algebra and
$(C,\delta)$ be a differential complex. If there exist two
morphisms of differential complexes
\[
\iota\colon (C[1],\delta_{[1]}) \to (V[1],q_1) \qquad
\text{and}
\qquad
\pi\colon (V[1],q_1)\to (C[1],\delta_{[1]})
\]
which are homotopy inverses, then there exist an $L_\infty$-algebra
structure $(C,\langle\,\rangle_1^{},\langle\,\rangle_2^{},\dots)$ on
$C$ extending its differential complex structure, and making
$(V,q_1,q_2,\dots)$ and
$(C,\langle\,\rangle_1^{},\allowbreak\langle\,\rangle_2^{},\dots)$ be
quasiisomorphic $L_\infty$-algebras via an
$L_\infty$-quasiisomorphism $\iota_\infty^{}$ extending $\iota$.
\end{theorem}
In case $\pi\imath=\operatorname{Id}_{C[1]}$, 
explicit formulas for such a transfer are described in
\cite{fuka,KonSoi} and \cite{Schuhmacher} in terms of
summation over rooted trees \cite[Definition 6]{KoSo}
\[ \langle\,\rangle^{}_n=\sum_{\Gamma\in
T_n}\varepsilon_\Gamma^{}\frac{Z_\Gamma(\imath,\pi,K,q_i)}{|\operatorname{Aut}\Gamma|},\]
where $K\in\Hom^{-1}(V[1],V[1])$ is an homotopy between $\iota\pi$ and 
$\operatorname{Id}_{V[1]}$,
$T_n$ is the set of rooted trees with $n$ tails, $\varepsilon_\Gamma^{}=\pm1$ is
a sign depending on the combinatorics of the tree $\Gamma$ and $
\operatorname{Aut}\Gamma$ is the group of automorphisms of $\Gamma$.
Each tail edge of a tree is
decorated by the operator
$\imath$, each internal edge is decorated by the
suspended operator
$K$ and the root edge is decorated by the
suspended operator
$\pi$; every internal vertex $v$ carries the  operation $q_r$,
where $r$ is the number of edges having $v$ as endpoint. Then $Z_{\Gamma}$ is the evaluation
of such a decorated graph according to the usual operadic
rules; see \cite[Thm.~2.3.1]{fuka} for an explicit recursive formula.

\bigskip
\section{The suspended mapping cone of $\chi\colon L\to M$.}
\label{sec.suspendedmap}

The suspended mapping cone of the DGLA morphism $\chi\colon
L\to M$ is the graded vector space
\[
C_\chi=\Cone(\chi)[-1],
\]
where $\Cone(\chi)=L[1]\oplus M$ is the mapping cone of
$\chi$. More explicitly,
\[ C_\chi=\mathop{\oplus}_iC_\chi^i,\qquad C_\chi^i=
L^i\oplus M^{i-1}.\]
The suspended mapping cone has  a natural differential
$\delta\in\Hom^1(C_{\chi},C_{\chi})$
given by
\[ \delta(l,m)=(dl,\chi(l)-dm),\qquad l\in L,m\in M.\]%

Denote by $\langle\;\rangle_1\in \Hom^1(C_{\chi}[1],C_{\chi}[1])$ the suspended differential,
namely 
\[ \langle (l,m)\rangle_1=(-dl,-\chi(l)+dm),\qquad l\in L,m\in M.\]%

\begin{remark}\label{remark.conequotient}
If $\chi$ is injective, then  the
projection on the second factor induces a quasiisomorphism of differential
complexes $\pi_2\colon C_\chi\to (M/\Image(\chi))[-1]$. In particular,  it
induces isomorphisms
$H^i(C_\chi)\simeq H^{i-1}(\coker(\chi))$, for every $i$.
\end{remark}

Setting
$M[t,dt]=\mathbb{K}[t,dt]\otimes M$, then
\[ H_\chi=\{(l,m(t,dt))\in L\times M[t,dt]\mid m(0,0)=0,\,
m(1,0)=\chi(l)\}\]%
is a differential graded Lie algebra.
The differential
on $H_\chi$ is $(l,m(t,dt))\mapsto(dl,dm(t,dt))$; since the
differential on $H_\chi$ has degree 1, the
suspended differential
$q_1\colon H_\chi[1]\to H_\chi[1]$ is the opposite
differential:
\[ {q}_1(l,m(t,dt))=-(dl,dm(t,dt)).
\]
    The integral operator $\int_a^b\colon
\mathbb{K}[t,dt]\to
\mathbb{K}$ extends naturally to a linear map of degree $-1$
\[ \int_a^b\colon M[t,dt]\to M,\qquad
\int_a^b(\sum_i t^im_i+t^idt\cdot n_i)=
\sum_i\left(\int_a^bt^idt\right)n_i.\]

\begin{lemma} The complexes $C_\chi[1]$ and $H_\chi[1]$ are homotopic; more precisely,
if one denotes by
\[ \imath\in\Hom^0(C_\chi[1],H_\chi[1]),\qquad \pi\in\Hom^0(
H_\chi[1],C_\chi[1]),\qquad K\in\Hom^{-1}(H_\chi[1],H_\chi[1])\]%
the linear maps defined as
\[  \imath(l,m)=(l,t\chi(l)+dt\cdot m),\qquad
\pi(l,m(t,dt))=\left(l,\int_0^1 m(t,dt)\right)
\]
\[
K(l,m)=\left(0,t\int_0^1m-\int_0^t m\right),
\]
then $\iota$ and
$\pi$ are morphisms of complexes and
\[ \pi\imath=\operatorname{Id}_{C_\chi[1]},\qquad
\operatorname{Id}_{H_\chi[1]}-\imath\,\pi=Kq_1+{q_1}K.
\]
\end{lemma}

\begin{proof} Straightforward.
\end{proof}

\bigskip
\section{The $L_\infty$ structure on $C_{\chi}$}
\label{sec:linfty}

By Quillen construction \cite{Qui},
the differential graded Lie algebra $H_\chi$
carries an $L_{\infty}$ structure
\[ {q}_k\colon \bigodot^k (H_\chi[1])\to H_\chi[1], \]
where ${q}_k=0$ for every $k\ge 3$,
\[ {q}_1(l,m(t,dt))=(-dl,-dm(t,dt))\]
and
\[q_2((l_1,m_1(t,dt))\odot(l_2,m_2(t,dt)
))=
(-1)^{\deg_{H_{\chi}}(l_1,m_1(t,dt))}([l_1,l_2],[m_1(t,dt),m_2(t,dt)]).\]

Results of Section \ref{sec.suspendedmap} tell us that we can apply
the homotopy transfer of structure theorem, to induce on $C_\chi$ an
$L_\infty$-algebra structure making $C_\chi$ and $H_\chi$ be
quasiisomorphic $L_\infty$-algebras.
Moreover, the linear maps of degree 1
\[ \langle\;\rangle_n^{}\colon \bigodot^n C_\chi[1]\to
C_\chi[1],\qquad n\ge 1,\]
defining the induced $L_\infty$-algebra structure on $C_\chi$ are
explicitly described in terms of
summation over rooted trees. In our case, the  properties
\[\pi^{} q_1 K=Kq_1\imath=0,\]
\[
q_2(\image{K}\otimes\image{K})\subseteq\ker\pi\cap\ker
K,
\qquad q_k=0\;\; \forall\; k\ge 3,\]
imply that, fixing the  number of tails, there exists at most one
isomorphism class of
trees  giving a nontrivial contribution.
\begin{itemize}

\item{} One tail: the only tree is
\[
\begin{xy}
,(-8,0)*{\circ};(0,0)*{\bullet}**\dir{-}?>*\dir{>}
,(0,0)*{\bullet};(8,0)*{\circ}**\dir{-}?>*\dir{>}
\end{xy}\qquad \rightsquigarrow\qquad
\begin{xy}
,(-13,0);(-8,0)*{\,\scriptstyle{\iota}\,}**\dir{-}
,(-8,0)*{\,\scriptstyle{\iota}\,};
(0,0)*{\,\scriptstyle{q_1}\,}**\dir{-}?>*\dir{>}
,(0,0)*{\,\scriptstyle{q_1}\,};
(8,0)*{\,\scriptstyle{\pi}\,}**\dir{-}
,(8,0)*{\,\scriptstyle{\pi}\,};
(16,0)**\dir{-}?>*\dir{>}
\end{xy}
\]
giving by operadic evaluation the formula
\[ \langle (l,m)\rangle_1=\pi
q_1\imath(l,m)=(-dl,-\chi(l)+dm).\]

\item{} Two tails:
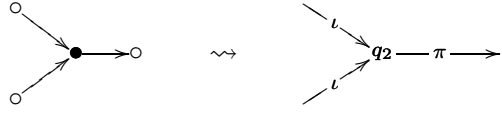
\[
\begin{xy}
,(-8,6)*{\circ};(0,0)*{\bullet}**\dir{-}?>*\dir{>}
,(-8,-6)*{\circ};(0,0)*{\bullet}**\dir{-}?>*\dir{>}
,(0,0)*{\bullet};(8,0)*{\circ}**\dir{-}?>*\dir{>}
\end{xy}
\qquad
\rightsquigarrow
\qquad
\begin{xy}
,(-10,6.66);(-6,4)*{\,\scriptstyle{\iota}\,}**\dir{-}
,(-10,-6.66);(-6,-4)*{\,\scriptstyle{\iota}\,}**\dir{-}
,(-6,4)*{\,\scriptstyle{\iota}\,};
(0,0)*{\,\,\scriptstyle{q_2}\,}**\dir{-}?>*\dir{>}
,(-6,-4)*{\,\scriptstyle{\iota}\,};(0,0)*{\,\,\scriptstyle{q_2}\,}**\dir{-}?>*\dir{>}
,(0,0)*{\,\,\scriptstyle{q_2}\,};
(8,0)*{\,\scriptstyle{\pi}\,}**\dir{-}
,(8,0)*{\,\scriptstyle{\pi}\,};
(16,0)**\dir{-}?>*\dir{>}
\end{xy}
\]
Again by operadic evaluation, this graph gives
\[ \langle
{\gamma_1}\odot{\gamma_2}\rangle_2=\pi
q_2(\imath({\gamma_1})
\odot\imath({\gamma_2})).\]

\item{} $n$ tails:
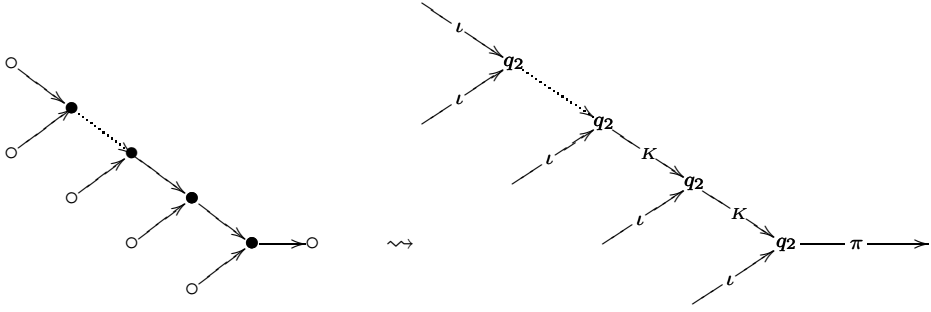
\[
\hskip -2.4em
\begin{xy}
,(-32,24)*{\circ};(-24,18)*{\bullet}**\dir{-}?>*\dir{>}
,(-32,12)*{\circ};(-24,18)**\dir{-}?>*\dir{>}
,(-24,6)*{\circ};(-16,12)*{\bullet}**\dir{-}?>*\dir{>}
,(-24,18)*{\circ};(-16,12)**\dir{.}?>*\dir{>}
,(-8,6)*{\circ};(0,0)*{\bullet}**\dir{-}?>*\dir{>}
,(-8,-6)*{\circ};(0,0)*{\bullet}**\dir{-}?>*\dir{>}
,(-16,12);(-8,6)*{\bullet}**\dir{-}?>*\dir{>}
,(-16,0)*{\circ};(-8,6)*{\bullet}**\dir{-}?>*\dir{>}
,(0,0)*{\bullet};(8,0)*{\circ}**\dir{-}?>*\dir{>}
\end{xy}
\qquad
\rightsquigarrow
\begin{xy}
,(-36,24)*{\,\,\scriptstyle{q_2}\,\,};(-24,16)*{\,\,\scriptstyle{q_2}\,\,}**\dir{.}?>*\dir{>}
,(-12,-8);(-7.2,-4.8)*{\,\scriptstyle{\iota}\,}**\dir{-}
,(-7.2,-4.8)*{\,\scriptstyle{\iota}\,};
(0,0)*{\,\,\scriptstyle{q_2}\,}**\dir{-}?>*\dir{>}
,(-24,0);(-19.2,3.2)*{\,\scriptstyle{\iota}\,}**\dir{-}
,(-19.2,3.2)*{\,\scriptstyle{\iota}\,};
(-12,8)*{\,\,\scriptstyle{q_2}\,\,}**\dir{-}?>*\dir{>}
,(-48,16);(-43.2,19.2)*{\,\scriptstyle{\iota}\,}**\dir{-}
,(-43.2,19.2)*{\,\scriptstyle{\iota}\,};
(-36,24)*{\,\,\scriptstyle{q_2}\,\,}**\dir{-}?>*\dir{>}
,(-48,32);(-43.2,28.8)*{\,\scriptstyle{\iota}\,}**\dir{-}
,(-43.2,28.8)*{\,\scriptstyle{\iota}\,};
(-36,24)*{\,\,\scriptstyle{q_2}\,\,}**\dir{-}?>*\dir{>}
,(-36,8);(-31.2,11.2)*{\,\scriptstyle{\iota}\,}**\dir{-}
,(-31.2,11.2)*{\,\scriptstyle{\iota}\,};
(-24,16)*{\,\,\scriptstyle{q_2}\,\,}**\dir{-}?>*\dir{>}
,(-12,8)*{\,\,\scriptstyle{q_2}\,\,};(-6,4)*{\,\scriptscriptstyle{K}\,}**\dir{-}
,(-6,4)*{\,\scriptscriptstyle{K}\,};
(0,0)*{\,\,\scriptstyle{q_2}\,}**\dir{-}?>*\dir{>}
,(-24,16)*{\,\,\scriptstyle{q_2}\,\,};(-18,12)*{\,\scriptscriptstyle{K}\,}**\dir{-}
,(-18,12)*{\,\scriptscriptstyle{K}\,};
(-12,8)*{\,\,\scriptstyle{q_2}\,\,}**\dir{-}?>*\dir{>}
,(0,0)*{\,\,\scriptstyle{q_2}\,};
(9.6,0)*{\,\scriptstyle{\pi}\,}**\dir{-}
,(9.6,0)*{\,\scriptstyle{\pi}\,};
(19.2,0)**\dir{-}?>*\dir{>}
\end{xy}
\]
\end{itemize}
This diagram gives, for every $n\ge 2$
the formula
\begin{multline*}
\langle{\gamma_1}\odot\cdots\odot
{\gamma_n}\rangle_n=\\
=\frac{(-1)^{n-2}}{2}\sum_{\sigma\in
S_n}\varepsilon(\sigma)\pi {q}_2(\imath(
{\gamma_{\sigma(1)}}
\!)\odot
Kq_2(\imath({\gamma_{\sigma(2)}})
\odot\cdots\odot
Kq_2(\imath({\gamma_{\sigma(n-1)}})
\odot\imath({\gamma_{\sigma(n)}}))\cdots)).
\end{multline*}
The
factor 1/2 in the above formula accounts for the cardinality
of the automorphism group of the graph involved.

\begin{remark}\label{rem.scalarextension}
The above construction of the $L_{\infty}$ structure on $C_{\chi}$
commutes
with tensor products of differential graded commutative
algebras.  This means that if $R$ is a DGCA, then the
$L_{\infty}$-algebra structure on the suspended mapping cone
of
$\chi\otimes{\rm id_R}\colon L\otimes R\to M\otimes R$ is naturally
isomorphic to the $L_{\infty}$-algebra
$C_{\chi}\otimes R$.
\end{remark}

A more refined description involving the original brackets in the
differential graded Lie algebras
$L$ and $M$ is obtained decomposing  the symmetric powers of
$C_\chi[1]$
into types:
\[
\bigodot^n \left(C_\chi[1]\right)=\bigodot^n
\Cone(\chi)=\bigoplus_{\lambda+\mu=n}
\left(\bigodot^\mu M\right)\otimes \left(\bigodot^\lambda
L[1]\right).
\]

The operation $\langle\,\rangle_2$ decomposes into
\begin{align*}
&{l_1}\otimes {l_2}\mapsto
(-1)^{\deg_L^{}(l_1)}[l_1,l_2]\in L;\qquad
{m_1}\otimes {m_2}\mapsto 0;\\
\\
&m\otimes l\mapsto
\dfrac{(-1)^{\deg_M^{}(m)+1}}{2}[m,\chi(l)]\in M
.
\end{align*}

For later use, we point out that, via decalage isomorphisms, the maps
$\langle-\rangle_1$ and $\langle-\rangle_2$ corresponds to

\[ \mu_1\in\Hom^1(C_{\chi},C_{\chi}),\qquad 
\mu_2\in\Hom^0(\bigwedge^2 C_{\chi}, C_{\chi}),\]

\[  \mu_1(l,m)=(dl,\chi(l)-dm),\]
\[  \mu_2((l_1,m_1)\wedge (l_2,m_2))=
\left([l_1,l_2],
\frac{1}{2}[m_1,\chi(l_2)]+\frac{(-1)^{\deg_L^{}(l_1)}}{2}[\chi(l_1),m_2]\right).\]

For every $n\geq 2$ it is easy to see that
$\langle{\gamma_1}\odot\cdots\odot{\gamma_{n+1}}
\rangle_{n+1}$
can be nonzero only
if the multivector
${\gamma_1}\odot\cdots\odot{\gamma_{n+1}}$
belongs to
$\bigodot^{n}M\otimes L[1]$.
For $n\ge 2$, $m_1,\ldots,m_{n}\in M$ and $l\in L[1]$ the
formula for
$\langle-\rangle_{n+1}$
described above becomes
\begin{multline*}
\langle {m_1}\odot\cdots\odot
{m_{n}}\odot l\rangle_{n+1}=\\
=(-1)^{n-1}\sum_{\sigma\in
S_{n}}\varepsilon(\sigma)\pi
q_2((dt) m_{\sigma(1)}\odot
Kq_2((dt)m_{\sigma(2)}\odot\cdots\odot
Kq_2((dt)m_{\sigma(n)}\odot
t\chi(l))\cdots)).
\end{multline*}
Define recursively
a sequence of polynomials $\phi_i(t)\in \Q[t]\subseteq \mathbb{K}[t]$
and rational numbers
$I_n$ by the rule
\[ \phi_1(t)=t,\qquad I_n=\int_0^1\phi_n(t)dt,\qquad
\phi_{n+1}(t)=\int_0^t\phi_{n}(s)ds-tI_n.\]
By the definition of the homotopy operator $K$ we have, for every
$m\in M$
\[ K((\phi_n(t)dt)m)=-\phi_{n+1}(t)m.\]
Therefore, for every  $m_1,m_2\in M$
   we have
\[
Kq_2((dt\cdot m_1)\odot
\phi_n(t){m_2})=(-1)^{\deg_M^{}(m_1)}\phi_{n+1}(t)
[m_1,m_2].\]
Therefore, we find:
\begin{multline*}
\langle m_1\odot\cdots\odot
m_{n}\odot
l\,\rangle_{n+1}=\qquad\qquad\qquad\\
\quad=(-1)^{n-1}\sum_{\sigma\in
S_{n}}\varepsilon(\sigma)\pi q_2((dt)
m_{\sigma(1)}\odot
Kq_2((dt)m_{\sigma(2)}
\odot\cdots\odot
Kq_2((dt)m_{\sigma(n)}\odot
t\chi(l))\cdots))\\
   =(-1)^{n-1+\deg_M(m_{\sigma(n)})}\sum_{\sigma\in
S_{n}}\varepsilon(\sigma)\pi q_2((dt)
m_{\sigma(1)}\odot
Kq_2((dt)m_{\sigma(2)}
\odot\cdots\odot
\phi_2(t)[m_{\sigma(n)},\chi(l)]\cdots))\\
\qquad =(-1)^{n-1+\sum_{i=2}^n\deg_M(m_{\sigma(i)})}\sum_{\sigma\in
S_{n}}\varepsilon(\sigma)\pi q_2((dt)
m_{\sigma(1)}\odot\phi_{n}(t)[m_{\sigma(2)},\cdots,
[m_{\sigma(n)},\chi(l)]\cdots])\\
=(-1)^{n+\sum_{i=1}^n\deg_M(m_i)}I_n\sum_{\sigma\in
S_{n}}\varepsilon(\sigma)
[m_{\sigma(1)},[m_{\sigma(2)},\cdots,
[m_{\sigma(n)},\chi(l)]\cdots]]\in M
\end{multline*}

\begin{theorem}\label{thm.coefficienti}
For every $n\ge 2$ we have
\begin{multline*}
\langle {m_1}\odot\cdots\odot
{m_{n}}\odot
l\,\rangle_{n+1}=\\
=-(-1)^{\sum_{i=1}^n\deg_M(m_i)}\frac{B_n}{n!}\sum_{\sigma\in
S_{n}}\varepsilon(\sigma)
[m_{\sigma(1)},[m_{\sigma(2)},\cdots,
[m_{\sigma(n)},\chi(l)]\cdots]],
\end{multline*}
where the $B_n$ are the Bernoulli numbers, defined by 
the series expansion identity
\[
\sum_{n=0}^{\infty}B_n\frac{x^n}{n!}=\frac{x}{e^x-1}
=1-\frac{x}{2}+\frac{x^2}{12}-\frac{x^4}{720}+\cdots
\]
\end{theorem}

\begin{proof} Since $B_{2k+1}=0$ for every $k>0$, 
it is sufficient to prove that $B_n=-n!I_n$ for every
$n\ge 1$.
Consider the polynomials $\psi_0(t)=1$ and
$\psi_n(t)=n!(\phi_n(t)-I_n)$ for $n\geq 1$. Then, for
any
$n\geq 1$,
\[
\frac{d}{dt}\psi_n(t)=n\cdot\psi_{n-1}(t),\qquad
\int_0^1\psi_n(t)dt=0.\]
Therefore  the $\psi_n(t)$ 
satisfy the recursive relations (see e.g. \cite{remmert})
of the
Bernoulli polynomials $B_n(t)$, 
defined by 
the series expansion identity
\[
\sum_{n=0}^{\infty}B_n(t)\frac{x^n}{n!}=\frac{xe^{tx}}{e^x-1}.
\]
In particular $B_n=B_n(0)=\psi_n(0)=-n!I_n$ for every $n\ge 1$.
\end{proof}

\begin{remark} 
Recently, the relevance of Bernoulli numbers in
deformation theory  has been also
remarked by Ziv Ran in \cite{ranATOMS}. In particular,  Ran's
``JacoBer" complex seems to
be closely related to the coderivation $Q$ defining the $L_\infty$
structure
on $C_\chi$.
\end{remark}

\bigskip
\section{The functors $\MC_\chi$ and $\Def_\chi$ revisited}
\label{sec.functorsrevisited}

Having introduced an $L_\infty$ structure on $C_\chi$ in
Section~\ref{sec:linfty}, we have a corresponding
Maurer-Cartan functor \cite{fuka,K}
$\MC_{C_\chi}\colon\mathbf{Art}\to\mathbf{Set}$,
defined as
\[\MC_{C_\chi}(A)=\left\{ \gamma\in
C_\chi[1]^0\otimes\mathfrak{m}_A
\;\strut\left\vert\;
\sum_{n\ge1}\frac{\langle\gamma^{\odot
n}\rangle_n}{n!}=0\right.\right\},\qquad
A\in \mathbf{Art}.
\]

Writing $\gamma=(l,m)$, with $l\in
L^1\otimes\mathfrak{m}_A$ and $m\in
M^0\otimes\mathfrak{m}_A$, the Maurer-Cartan equation becomes
\begin{align*}
0&=\sum_{n=1}^\infty\frac{\langle(l,m)^{\odot
n}\rangle_n}{n!}\\
&=
\langle (l,m)\rangle_1+\frac{1}{2}\langle
l^{\odot 2}\rangle_2+\langle
m\otimes l\rangle_2+\frac{1}{2}\langle
m^{\odot 2}\rangle_2+\sum_{n\geq 2}
\frac{n+1}{(n+1)!}\langle m^{\odot n}\otimes
l\rangle_{n+1}\\
&=
\left(-dl-\frac{1}{2}[l,l],
-\chi(l)+dm -\frac{1}{2}[m,\chi(l)]
+\sum_{n\geq 2}
\frac{1}{n!}\langle m^{\odot n}\otimes
l\rangle_{n+1}\right)\in(L^2\oplus M^1)\otimes\mathfrak{m}_A.
\end{align*}
According to Theorem~\ref{thm.coefficienti}, 
since $\deg_M(m)=\deg_{C_{\chi[1]}}(m)=0$, we have 
\[ \langle m^{\odot n}\otimes
l\rangle_{n+1}=-\frac{B_n}{n!}
\sum_{\sigma\in
S_{n}}[m,[m,\cdots,
[m,\chi(l)]\cdots]]=-B_n\ad_m^n(\chi(l)).\]
The Maurer-Cartan equation on $C_{\chi}$ 
is therefore  equivalent
to
\[
\begin{cases}
\displaystyle{dl+\frac{1}{2}[l,l]=0}\\
\\
\displaystyle{\chi(l)-dm+\frac{1}{2}[m,\chi(l)]+
\sum_{n=2}^{\infty}\frac{B_n}{n!} \ad_m^n(\chi(l))=0}.
\end{cases}
\]
Since $B_0=1$ and $B_1=-\dfrac{1}{2}$, we can write the
second equation as
\[
[m,\chi(l)]-dm+
\sum_{n=0}^{\infty}\frac{B_n}{n!}\ad_m^n(\chi(l))
=[m,\chi(l)]-dm+\frac{\ad_m}{e^{\ad_m}-1}(\chi(l))=0.
\]
Applying the invertible operator $\dfrac{e^{\ad_m}-1}{\ad_m}$  we get
\[
0=\chi(l)+\dfrac{e^{\ad_m}-1}{\ad_m}([m,\chi(l)]-dm)=e^m\ast\chi(l).\]

Therefore, the
Maurer-Cartan equation for the $L_\infty$-algebra
structure on $C_\chi$ is equivalent to
\[
\begin{cases}
dl+\dfrac{1}{2}[l,l]=0\\
e^m\ast\chi(l)=0
\end{cases}
\]
and the Maurer-Cartan functor $\MC_\chi$ described in the
introduction
is precisely the Maurer-Cartan
functor corresponding to the
$L_\infty$ structure on $C\chi$.
\medskip

Recall that
the deformation functor
associated to an $L_\infty$-algebra
$\mathfrak{g}$ is
$\Def_{\mathfrak{g}}=\MC_{\mathfrak{g}}/\sim$, where
$\sim$ denotes homotopy equivalence of solutions of the
Maurer-Cartan equation: two elements
$\gamma_0$ and
$\gamma_1$ of
$\MC_{\mathfrak g}(A)$ are called homotopy equivalent if
there exists an element $\gamma(t,dt)\in\MC_{{\mathfrak
g}[t,dt]}(A)$ with
$\gamma(0)=\gamma_0$ and
$\gamma(1)=\gamma_1$.

We have already  identified the functor
$\MC_{C_\chi}$ with the functor
$\MC_\chi$. Now we want  to show that, under this identification, the
homotopy equivalence on
$\MC_{C_\chi}$ is the same thing of gauge equivalence on
$\MC_\chi$ described in the introduction, so that
\[
\Def_\chi\simeq\Def_{C_\chi}.
\]
We will need the
following lemma (see the Appendix A for a proof).

\begin{lemma}\label{lemma:e-star-g}
Let $\mathfrak{g}$ be a differential graded Lie algebra and
let $A\in\mathbf{Art}$. Then, for any $x$ in
$\MC_{\mathfrak{g}}(A)$ and any
$g(t)\in \mathfrak{g}^0[t]\otimes\mathfrak{m}_A$, with
$g(0)=0$, the element $e^{g(t)}\ast x$
is an element of
$\MC_{{\mathfrak g}[t,dt]}(A)$.
Moreover all the elements of
$\MC_{\mathfrak{g}[t,dt]}(A)$ are obtained in this way.
\end{lemma}

We first show that homotopy implies gauge. Let
$(l_0,m_0)$ and $(l_1,m_1)$ be homotopy equivalent
elements of $\MC_{C_\chi}(A)$. Then there exists an element
$(\tilde{l},\tilde{m})$  of
$\MC_{C_\chi[s,ds]}(A)$ with
$(\tilde{l}(0),\tilde{m}(0))=(l_0,m_0)$ and
$(\tilde{l}(1),\tilde{m}(1))=(l_1,m_1)$.
According to Remark~\ref{rem.scalarextension},
the Maurer-Cartan equation for $(\tilde{l},\tilde{m})$ is
\[
\begin{cases}
d\tilde{l}+\dfrac{1}{2}[\tilde{l},\tilde{l}]=0\\
e^{\tilde{m}}\ast\chi(\tilde{l})=0
\end{cases}
\]
The first of the two equations above tells us that
$\tilde{l}$ is a solution of the Maurer-Cartan equation for
$L[s,ds]$. So, by Lemma~\ref{lemma:e-star-g}, there exists a
degree zero element
$\lambda(s)$ in $L[s]\otimes\mathfrak{m}_A$ with
$\lambda(0)=0$ such that $\tilde{l}=e^{\lambda}\ast l_0$.
Evaluating at $s=1$ we find $l_1=e^{\lambda_1}\ast l_0$. As a
consequence of $\tilde{l}=e^{\lambda}\ast l_0$, we also have
$\chi(\tilde{l})=e^{\chi(\lambda)}\ast\chi(l_0)$. Set
$\tilde{\mu}=\tilde{m}\bullet\chi(\lambda)\bullet m_0$, so
that
$\tilde{m}=\tilde{\mu}\bullet m_0\bullet(-\chi(\lambda))$ and
the second Maurer-Cartan equation is reduced to
$e^{\tilde{\mu}}\ast(e^{m_0}\ast\chi(l_0))=0$,
i.e., to
$e^{\tilde{\mu}}\ast 0=0$,
where we have used the fact that $(l_0,m_0)$ is a solution of
the Maurer-Cartan equation in $C_\chi$. This last equation is
equivalent to the equation $d\tilde{\mu}=0$ in
$(C_\chi[s,ds])^0\otimes\mathfrak{m}_A$. If we write
$\tilde{\mu}(s,ds)=\mu^0(s)+ds\,\mu^{-1}(s)$, then the
equation
$d\tilde{\mu}=0$ becomes
\[
\begin{cases}
\dot{\mu}^0-d^{}_M\mu^{-1}=0\\
d^{}_M\mu^0=0,
\end{cases}
\]
where $d^{}_M$ is the differential in the DGLA $M$.
The solution is, for any fixed $\mu^{-1}$,
\[
\mu^0(s)=\int_0^s d\sigma\, d^{}_M\mu^{-1}(\sigma)=-d^{}_M
\int_0^s d\sigma\,
\mu^{-1}(\sigma)
\]
Set $\nu=-\int_0^1 d s\,
\mu^{-1}(s)$. Then $m_1=\tilde{m}(1)=(d_M\nu)\bullet
m_0\bullet (-\chi(\lambda_1))$. Summing up, if $(l_0,m_0)$
and
$(m_1,l_1)$ are homotopy equivalent, then there exists
$(d\nu,\lambda_1)\in (d M^{-1}\otimes\mathfrak{m}_A)\times
(L^0\otimes\mathfrak{m}_A)$ such that
\[
\begin{cases}
l_1=e^{\lambda_1}\ast l_0\\
m_1=d\nu\bullet m_0\bullet (-\chi(\lambda_1)),
\end{cases}
\]
i.e., $(l_0,m_0)$
and
$(m_1,l_1)$ are gauge equivalent.

We now show that gauge
implies homotopy. Assume  $(l_0,m_0)$
and
$(m_1,l_1)$ are gauge equivalent. Then then there exist
$(d\nu,\lambda_1)\in (d M^{-1}\otimes{\mathfrak m})\times
(L^0\otimes{\mathfrak m})$ such that
\[
\begin{cases}
l_1=e^{\lambda_1}*l_0\\
m_1=d\nu\bullet m_0\bullet (-\chi(\lambda_1)).
\end{cases}
\]
Set $\tilde{l}(s,ds)=e^{s\lambda_1}*l_0$. By
Lemma~\ref{lemma:e-star-g},
$\tilde{l}$ satisfies
the equation $d\tilde{l}+\frac{1}{2}[\tilde{l},\tilde{l}]=0$.
Set $\tilde{m}=(d(s\nu))\bullet
m_0\bullet(-\chi(s\lambda_1))$. Reasoning as above, we find
\[
e^{\tilde m}*\chi(\tilde{l})=e^{d(s\nu)}*0=0.
\]
Therefore, $(\tilde{l},\tilde{m})$ is a solution of the
Maurer-Cartan equation in $C_\chi[s,ds]$. Moreover
$\tilde{l}(0)=l_0$, $\tilde{l}(1)=l_1$, $\tilde{m}(0)=m_0$
and
$\tilde{m}(1)=d\nu\bullet m_0\bullet (-\chi(\lambda_1))=m_1$,
i.e. $(l_0,m_0)$
and
$(m_1,l_1)$ are homotopy equivalent.
\\
Summing up, we have
\[
\Def_{C_\chi}=\frac{\MC_{C_\chi}}{\text{homotopy}}\simeq
\frac{\MC_{\chi}}{\text{gauge}}=\Def_{\chi}.
\]

\bigskip
\section{Functoriality}\label{sec:functoriality}

In the above section we have shown how to a morphism of
differential graded Lie algebras $\chi\colon L\to M$ is
associated a canonical  $L_\infty$ structure on
$C_\chi$. We will now discuss the functorial aspects of this
construction. Denote by $\mathbf{L}_\infty$ the category of
$L_\infty$-algebras and by $\mathbf{M}$ the
category of morphisms of differential graded Lie algebras; objects
in $\mathbf{M}$ are DGLA morphisms $\chi\colon L\to M$;
morphisms in $\mathbf{M}$ are commutative squares
\[
\xymatrix{
      L_1 \ar[r]^{f_L}
\ar[d]_{\chi^{}_1} &
L_2 \ar[d]^{\chi^{}_2}\\
    M_1  \ar[r]^{f_M} &
M_2\\
    }
\]
of DGLA morphisms.\\

It is immediate to observe that the above commutative
square induces a linear $L_{\infty}$-morphism $(f_L,f_M)\colon
C_{\chi^{}_1}\to
C_{\chi^{}_2}$ and then
$C\colon\mathbf{M}\to\mathbf{L}_\infty$ is a functor.
Moreover both $\MC$ and $\Def$ are functors from the category 
$\mathbf{L}_\infty$ to the category
of functors of Artin rings \cite{K,EDF,defomanifolds}. The functor 
$\MC$ acts on the morphisms of
$\mathbf{L}_\infty$ in the following way:
let $f_{\infty}\colon V\to W$ be an $L_{\infty}$-morphism, then
\[ \MC_{f_{\infty}}\colon \MC_{V}\to \MC_W \]
is the natural transformation given, for $A\in\mathbf{Art}$ and $v\in 
\MC_V(A)\subseteq
V[1]^0\otimes\mathfrak{m}_A$,  by
\[\MC_{f_{\infty}}(v)=
\sum_{n=1}^{\infty}\frac{1}{n!}f_{n}(v^{\odot n})\in
W[1]^{0}\otimes \mathfrak{m}_{A}.\]
The natural transformation $\MC_{f_{\infty}}$ preserves the homotopy 
equivalence and then
induces a natural transformation
\[ \Def_{f_{\infty}}\colon \Def_{V}\to \Def_W.\]

Recall from Section \ref{sec.linfty.intro} that an
$L_\infty$-morphism $f^{}_\infty$ is called a quasiisomorphism of
$L_\infty$-algebras if its linear part $f_1$ is a quasiisomorphism of
differential
complexes. Using the fact that every quasiisomorphism  of
$L_{\infty}$-algebras
induces an isomorphism of the associated  deformation functors
\cite{K}, the next
theorem becomes evident.

\begin{theorem}[\cite{semireg}]\label{thm.basic}
Consider a commutative diagram
of morphisms of differential
graded Lie algebras
\[
\xymatrix{
      L_1 \ar[r]^{f_L}
\ar[d]_{\chi^{}_1} &
L_2 \ar[d]^{\chi^{}_2}\\
    M_1  \ar[r]^{f_M} &
M_2\\
    }
\]
and assume that
$(f_L,f_M)\colon C_{\chi^{}_1}\to C_{\chi^{}_2}$ is a
quasiisomorphism of complexes (e.g. if both $f_L$ and $f_M$ are
quasiisomorphisms).
Then the natural transformation $\Def_{\chi^{}_1}\to
\Def_{\chi^{}_2}$ is an isomorphism.
\end{theorem}

\bigskip
\section{Cartan homotopies}
\label{sec.cartan}

In this section we formalize, under the notion of \emph{Cartan homotopy},
a set of standard identities that often arise in algebra and geometry
\cite[Appendix B]{clemens}.

\begin{definition}\label{def:cartan}
Let $(L,d,[,])$ and $(M,d,[,])$ be two differential graded Lie
algebras and
denote by $\delta$ the standard differential on $\Hom^*(L,M)$.
A linear map
$i\in \Hom^{-1}(L,M)$ is called a \emph{Cartan homotopy} if for every
$a,b\in L$
we have:
\[i([a,b])=[i(a),\delta i(b)],\qquad [i(a),i(b)]=0.\]
\end{definition}

Notice that, according to the definition of $\delta$, for every $a\in
L$ we have
\[\delta i(a)=d(i(a))+i(da).\]

For later use we point out that
$[i(a),[i(b),\delta i(c)]]=[i(a),i([b,c])]=0$ for every
$a,b,c$. It is moreover easy to verify that $\delta i$ is a morphism
of differential graded Lie algebras  and

\[i([a,b])=(-1)^{\deg(a)}[\delta i(a),i(b)]=\frac{1}{2}[i(a),\delta
i(b)]+
\frac{(-1)^{\deg(a)}}{2}[\delta i(a),i(b)].\]

\begin{example}\label{exa.realcartan}
The name Cartan homotopy has a clear origin in differential geometry.
Namely, let $M$
be a differential manifold, ${\mathcal X}(M)$ be the Lie algebra of
vector fields on
$M$, and ${\mathcal E}nd^*(\Omega^*(M))$ be the Lie algebra of
endomorphisms of the
de Rham algebra of $M$. The Lie algebra ${\mathcal X}(M)$ can be seen
as a DGLA
concentrated in degree zero, and the graded Lie algebra ${\mathcal
E}nd^*(\Omega^*(M))$ has a degree one differential given by
$[d,-]$, where $d$
is the de Rham differential. Then the contraction
\[
i\colon {\mathcal X}(M)\to {\mathcal E}nd^*(\Omega^*(M))[-1]
\]
is a Cartan homotopy and its differential is the Lie derivative
\[
\delta i={\mathcal L}\colon {\mathcal X}(M)\to {\mathcal
E}nd^*(\Omega^*(M)).
\]
In fact, by classical Cartan's homotopy formulas
\cite[Section 2.4]{AbrahamMarsden}, for any two vector fields $X$ and
$Y$ on $M$, we
have
\begin{enumerate}
\item ${\mathcal L}_X=di_X+i_Xd=[d,i_X]$;

\item $i_{[X,Y]}={\mathcal L}_X i_Y-i_Y{\mathcal L}_X=[{\mathcal
L}_X,i_Y]=[i_X,{\mathcal L}_Y]$;

\item $[i_X,i_Y]=0$.
\end{enumerate}
Note that the first Cartan formula above actually states that
$\delta i={\mathcal L}$. Indeed  ${\mathcal X}(M)$ is
concentrated in  degree
zero and then its differential is trivial.
\end{example}

\begin{example}\label{ex.cambiobasepercartan}
The composition of a Cartan homotopy with a morphism of DGLAs is a 
Cartan homotopy.
If $i\colon L\to M[-1]$ is a Cartan
homotopy and
$\Omega$ is a differential graded-commutative algebra, then
its natural extension
\[i\otimes \operatorname{Id}\colon L\otimes \Omega\to (M\otimes \Omega)[-1],\qquad
a\otimes \omega\mapsto i(a)\otimes\omega,\]
is  a Cartan homotopy.
\end{example}

\begin{proposition}\label{prop.cartan}
Let $i\colon L\to M[-1]$ be a Cartan homotopy and  $\chi=\delta
i\colon L\to M$.
Then the linear map
\[ \tilde{i}\colon L\to C_{\chi},\qquad \tilde{i}(a)=(a,i(a))\]
is a linear  $L_{\infty}$-morphism.
\end{proposition}

\begin{proof} By  decalage isomorphism, the $L_{\infty}$ structure on
$C_{\chi}$ is given by the higher brackets  
$\mu_n\in\Hom^{2-n}(\bigwedge^nC_{\chi},C_{\chi})$, $n\ge 1$, where
\[ \mu_1((l,m))=(dl,\chi(l)-dm),\]
\[ \mu_2((l,m)\wedge (h,k))=
\left([l,h],\frac{1}{2}[m,\chi(h)]+\frac{(-1)^{\deg(l)}}{2}[\chi(l),k]\right)\]

and for $n\ge 3$
\[
\mu_n((l_1,m_1)\wedge\cdots\wedge(l_n,m_n))=\frac{B_{n-1}}{(n-1)!}\sum_{\sigma\in
S_n}\pm
[m_{\sigma(1)},[\cdots,[m_{\sigma(n-1)},\chi(l)_{\sigma(n)}]\cdots]]
\]
It is straightforward to check that $\tilde{i}$
commutes with every bracket, i.e.
\[ \tilde{i}(dx)=\mu_1(\tilde{i}(x)),\qquad
\tilde{i}([x,y])=\mu_2(\tilde{i}(x)\wedge
\tilde{i}(y)),\] and for  $n\ge 3$
\[ \mu_n(\tilde{i}(x_1)\wedge\cdots\wedge \tilde{i}(x_n))=0.\]
Therefore $\tilde{i}$ is a linear $L_{\infty}$-morphism.
\end{proof}

\bigskip
\section*{Appendix A: gauge vs. homotopy}

In this Appendix we briefly discuss the relation
between homotopy and gauge equivalence for solutions of the
Maurer-Cartan equation for a given differential graded Lie
algebra $L$. We also give a proof of Lemma~\ref{lemma:e-star-g},
which is here presented as
Corollary A.1.

\begin{propositiona1}\label{prop.hvg}
Let $(L,d,[~,~])$ be a differential graded Lie algebra such that:
\begin{enumerate}

\item $L=M\oplus C\oplus D$ as graded vector spaces.

\item $M$ is a differential graded subalgebra of $L$.

\item  $d\colon C\to D[1]$ is an isomorphism of graded
vector spaces.
\end{enumerate}
Then, for every
$A\in \mathbf{Art}$ there exists a bijection
\[ \alpha\colon \MC_M(A)\times (C^0\otimes
\mathfrak{m}_A)\mapor{\sim}\MC_L(A),\qquad (x,c)\mapsto e^c\ast
x.\]
\end{propositiona1}

\begin{proof} This is essentially proved in \cite[Section 5]{SchSta}
using
induction on the length of $A$ and the Baker-Campbell-Hausdorff
formula.\\ Here we sketch a simpler  proof based on formal
theory of deformation functors \cite{Sch,Rim,FM1,ManettiDGLA}.\\
The map $\alpha$ is a natural
transformation of homogeneous functors,
so it is sufficient to show that $\alpha$
is bijective on tangent spaces and injective on obstruction spaces.
Recall that the tangent space of $\MC_L$ is $Z^1(L)$, while its
obstruction space is contained in $H^2(L)$.
The functor $A\mapsto C^0\otimes\mathfrak{m}_A$
is smooth with tangent space $C^0$ and therefore tangent and
obstruction spaces of the functor
\[ A\mapsto \MC_M(A)\times(C^0\otimes\mathfrak{m}_A)\]
are respectively $Z^1(M)\oplus C^0$ and $H^2(M)$.
The tangent map is
\[ Z^1(M)\oplus C^0\ni (x,c)\mapsto e^c\ast x=x-dc\in
Z^1(M)\oplus d(C^0)=Z^1(M)\oplus D^1=Z^1(L)\]%
and it is an isomorphism.
The inclusion $M\hookrightarrow L$ is a quasiisomorphism,
therefore the obstruction to lifting $x$ in $M$ is equal to
the obstruction to lifting $x=e^0\ast x$ in $L$. We conclude
the proof by observing that, according to
\cite[Prop. 7.5]{FM1}, \cite[Lemma 2.21]{ManettiDGLA},
the obstruction maps of
Maurer-Cartan functor are invariant under the gauge action.
\end{proof}

\begin{corollarya1}\label{cor:homvsgauge1}
Let $M$ be a differential graded Lie algebra, $L=M[t,dt]$ and
$C\subseteq M[t]$ the
subspace consisting of polynomials $g(t)$ with $g(0)=0$. Then for
every
$A\in \mathbf{Art}$
the map
$(x,g[t])\mapsto e^{g(t)}\ast x$ induces an isomorphism
\[\MC_M(A)\times (C^0\otimes \mathfrak{m}_A)\simeq\MC_L(A).\]
\end{corollarya1}
\begin{proof} The data $M,C$ and $D=d(C)$ satisfy the condition of
Proposition~A.1.\end{proof}

\begin{corollarya2} Let $M$ be a differential graded Lie algebra.
Two elements
$x_0,x_1\in \MC_M(A)$ are gauge equivalent if and
only if they are homotopy equivalent.
\end{corollarya2}
\begin{proof}
If $x_0$ and $x_1$ are gauge equivalent, then
there exists $g\in M^0\otimes\mathfrak{m}_A$ such that
$e^g\ast x_0=x_1$. Then, by Corollary
A.1. $x(t)=e^{t\, g}\ast x_0$ is an element of
$\MC_{M[t,dt]}(A)$ with $x(0)=x_0$ and $x(1)=x_1$, i.e.,
$x_0$ and $x_1$ are homotopy equivalent.\\
Vice versa, if $x_0$ and $x_1$ are homotopy equivalent, there exists
$x(t)\in \MC_{M[t,dt]}(A)$ such that $x(0)=x_0$ and
$x(1)=x_1$. By Corollary A.1., there exists
$g(t)\in M^0[t]\otimes \mathfrak{m}_A$ with $g(0)=0$ such that
$x(t)=e^{g(t)}\ast x_0$. Then $x_1=e^{g(1)}\ast x_0$, i.e., $x_0$
and $x_1$ are gauge equivalent.
\end{proof}

\bigskip
\part{dg-Grassmann functors}

Let
$W$ be vector space over $\K$. The total Grassmannian
of $W$ is
\[
\Grass(W)=\{\text{linear subspaces of }W\}.
\]
The group $\Aut(W)$ of linear automorphisms of $W$ acts on
$\Grass(W)$. Denoting by $\Grass(V,W)\subseteq \Grass(W)$
the orbit  of a subspace $V\subseteq W$
we have
\[
\Grass(V,W)=\frac{\Aut(W)}{\Aut(V,W)},\quad
\text{where}\quad
\Aut(V,W)=\{g\in \Aut(W) \mid g(V)= V\}.
\]

The infinitesimal neighborhood of $V$ in $\Grass(W)$
is the formal moduli space for the functor
$\Grass_{V,W}\colon \mathbf{Art}\to\mathbf{Set}$:
\begin{align*}
\Grass_{V,W}(A)=&\{\phi\colon {\rm Spec}(A)\to \Grass(W)\mid
\phi({\rm Spec}(\mathbb{K}))=V\}\\
=&\{\text{free $A$ submodules }V_A\subseteq W\otimes A\mid V_A\otimes_A\K=V\}\\
=&\{ f(V\otimes A)\subseteq W\otimes A\mid f\in \Aut(W\otimes A),\;
f_{|W}=\operatorname{Id}\}.
\end{align*}

Since $\{f\in \Aut(W\otimes A)\mid f_{|W}=\operatorname{Id}\}=
\exp(\mathfrak{gl}(W)\otimes\mathfrak{m}_A)$ we can write
\[
\Grass_{V,W}(A)=\frac{\exp(\mathfrak{gl}(W)\otimes
\mathfrak{m}_A)}{\exp(L^0_{V,W}\otimes \mathfrak{m}_A)}
\]
where
\[
L^0_{V,W}=\{g\in \mathfrak{gl}(W)\mid g(V)\subseteq V\}\]
and the action is given by
\[
\exp(L^0_{V,W}\otimes \mathfrak{m}_A)\times
\exp(\mathfrak{gl}(W)\otimes
\mathfrak{m}_A)
\to \exp(\mathfrak{gl}(W)\otimes
\mathfrak{m}_A),\qquad
(e^a,e^m)\mapsto e^me^{-a}.
\]
In conclusion, the functor $\Grass_{V,W}$ coincides with the
deformation functor
$\Def_{\chi}$ where $\chi\colon L^0_{V,W}\to \mathfrak{gl}(W)$
is the inclusion. Indeed the differential graded Lie algebras
$L^0_{V,W}$ and $\mathfrak{gl}(W)$ are
concentrated in degree 0,
   $\MC_{\chi}(A)=\exp(\mathfrak{gl}(W)\otimes
\mathfrak{m}_A)$ and the gauge action is given by
\[ \exp(L^0_{V,W}\otimes \mathfrak{m}_A)\times
\exp(\mathfrak{gl}(W)\otimes \mathfrak{m}_A)
\to \exp(\mathfrak{gl}(W)\otimes \mathfrak{m}_A),\qquad
(e^a,e^m)\mapsto e^me^{-\chi(a)}.\]

\bigskip
\section{The coarse  dg-Grassmannian}
\label{sec.coarsegrass}

The considerations of the above section suggest the following
generalization from vector spaces to differential complexes.
Let
$(W,d)$ be a dg-vector space and denote by:
\begin{enumerate}

\item $\Aut(W)$ the group of automorphisms of the graded vector space $W$;

\item $\Aut(W,d)$ the group of automorphisms of the differential
graded
vector space $(W,d)$, i.e. the subgroup of $\Aut(W,d)$ consisting of
linear
automorphisms which commute with the differential $d$;

\item $\Aut^0(W,d)$ the subgroup of $\Aut(W,d)$ of automorphisms
inducing the identity in cohomology.

\end{enumerate}

Define the coarse Grassmannian $\Grass(W)$ as the quotient
\[ \Grass(W)=\frac{G(W)}{\Aut^0(W,d)},\qquad \text{where}\qquad
G(W)=\{\text{subcomplexes of }(W,d)\,\}.\]%
Note that $\Aut(W,d)$ acts on $G(W)$ and then the quotient group
$\Aut(W,d)/\Aut^0(W,d)$ acts on $\Grass(W)$.\\

If $d=0$ then $\Aut^0(W,d)=\{\operatorname{Id}\}$ and then $\Grass(W)$
is the standard Grassmannian. Notice that the cohomology functor gives a
map
\[ h\colon \Grass(W)\to \Grass(H^*(W)),\qquad h(V)=\image(H^*(V)\to H^*(W)).\]

Denoting by $\Grass(W)^s\subseteq \Grass(W)$ the ``open'' subset consisting of
subcomplexes $V\subseteq W$ such that
$H^*(V)\to H^*(W)$ is injective, we shall prove later
(Theorems~\ref{thm.pregrassinjectivecohomology} and 
\ref{thm.grassinjectivecohomology})
that, is some sense,
the map $h\colon \Grass(W)^s\to \Grass(H^*(W))$ is a local isomorphism;
this fact justifies the quotient of $G(W)$ by the action of 
$\Aut^0(W,d)$.\\

Given a subcomplex $V\subseteq W$ we denote $\Grass(V,W)=G(V,W)/\Aut^0(W)$,
where $G(V,W)$ the set of subcomplexes of
$W$ that are isomorphic to $V$ as graded vector spaces; equivalently 
$U\in G(V,W)$ if and only if
$U\in G(W)$ and there exists $f\in \Aut(W)$ such that $f(V)=U$ and 
then we have a natural
identification
\[ \frac{\{f\in \Aut(W)\mid df(V)\subseteq f(V)\}}{\Aut(V,W)}\;\mapor{\simeq}\;
G(V,W),\qquad f\mapsto f(V),\]
where $\Aut(V,W)$ is the subgroup of $\Aut(W)$ consisting of the
automorphisms $g$ such that $g(V)=V$.

\begin{remark}
To put some structure on $\Grass(W)$
the natural choice is to take
$\Grass(V,W)$ as its components; then  try to
define $\Grass(V,W)$ as  categorical quotient
(in a suitable category: varieties, schemes,
stacks, ...) using the identification
\[ \Grass(V,W)=\frac{\{f\in \Aut(W)\mid df(V)\subseteq 
f(V)\}}{\Aut^0(W,d)\times \Aut(V,W)}.\]
Note that $(\varphi,\psi)\in
\Aut^0(W,d)\times \Aut(V,W)$ acts of $\{f\in \Aut(W)\mid 
df(V)\subseteq f(V)\}$ by the formula
$(\varphi,\psi)\cdot f=\varphi f \psi^{-1}$.\\
Unfortunately we may not expect that $\Grass(V,W)$
is separated in general.
Consider $W^0=W^1=\K\oplus \K$, $W^i=0$ for $i\not=0,1$, and
$d\colon W^0\to W^1$ of rank 1.\\
If $V^0=V^1=\K$, then $V$ is a subcomplex of $W$ if and only if
$V^0=\ker(d)$ or $V^1=\image(d)$.
Therefore $G(V,W)\subseteq \Proj(W^0)\times
\Proj(W^1)=\Proj^1\times\Proj^1$ is the union of
two intersecting lines. In particular  $G(V,W)$ is singular at
the point $V^0=\ker(d)$,
$V^1=\image(d)$.\\
The group  $\Aut^0(W,d)$ acts transitively on $\{V\mid
V_0\not=\ker(d)\}$,
$\{V\mid V_1\not=\image(d)\}$ and then $\Grass(V,W)$
contains three point and it is not separated (same type of
$\{xy=0\}/\K^*$).
\end{remark}

\bigskip

\section{Infinitesimal study}

Let $(W,d)$ be a complex of   vector spaces and $V\subseteq W$ a subcomplex;
denote by $L_{W}$ and $L_{V,W}$
the differential graded Lie algebras
\[ L_W=\Hom^*(W,W),\qquad L_{V,W}=\{g\in\Hom^*(W,W)\mid g(V)\subseteq V\}.\]%

An Artinian algebra can be seen as a
differential
complex with trivial differential; for every
$A\in\mathbf{Art}$ we still denote by $d$ the differential on
$W\otimes A$.
Since the differential on $A$ is trivial, we have
\[d\colon W^i\otimes A\to W^{i+1}\otimes A,\qquad d(v\otimes
a)=d(v)\otimes a.\]

Let $V\subseteq W$ be a subcomplex and consider the functors
\[\Aut_W,\Aut_{W,d},\Aut^0_{W,d}, \Aut_{V,W}\colon\mathbf{Art}\to
\mathbf{Groups}\]
defined as
\begin{align*}
\Aut_W(A)&=\{f\in \Hom^0_A(W\otimes A,W\otimes A)\mid
f\equiv  \operatorname{Id}\pmod{\mathfrak{m}_A}\}\\
\Aut_{W,d}(A)&=\{f\in \Aut_{W}(A)\mid fd=df\}\\
    \Aut_{V,W}(A)&=\{f\in
\Aut_W(A)\mid f(V\otimes A)=V\otimes A\}\\
    \Aut^0_{W,d}(A)&=\{f\in
\Aut_{W,d}(A)\mid H^*(f) \text{ is the identity  on  } H^*(W\otimes
A,d)\}
\end{align*}

The above functors  are smooth and homogeneous \cite{Rim,ManettiDGLA}. Their
tangent
spaces are
\begin{align*}
T^1\Aut_W&=\Hom^0(W,W)=L^0_W\\
T^1\Aut_{W,d}&=\{f\in \Hom^0(W,W)\mid fd=df\}=Z^0(\Hom^*(W,W))\\
T^1\Aut_{W,V}&=\{f\in \Hom^0(W,W)\mid  f(V)\subseteq
V\}=L_{V,W}^0\\
T^1\Aut^0_{W,d}&=\{f\in \Hom^0(W,W)\mid fd=df,\; f(\ker(d))\subseteq
\image(d)\}
\end{align*}

\begin{lemma}\label{lem.expo}
For every $A\in\mathbf{Art}$, the exponential map gives isomorphisms
\begin{align*} \exp\colon& \Hom^0(W,W)\otimes\mathfrak{m}_A\to
\Aut_W(A)\\
\exp\colon& Z^0(\Hom^*(W,W))\otimes\mathfrak{m}_A\to \Aut_{W,d}(A)\\
\exp\colon& B^0(\Hom^*(W,W))\otimes\mathfrak{m}_A\to \Aut^0_{W,d}(A)\\
\exp\colon& L^0_{V,W}\otimes\mathfrak{m}_A\to \Aut_{V,W}(A)
\end{align*}

\end{lemma}
\begin{proof}
First we note that
For every $i$ there exists a natural exact sequence
\[ 0\to B^i(\Hom^*(W,W))\to Z^i(\Hom^*(W,W))\to
\Hom^i(H^*(W),H^*(W))\to 0\]
and therefore
\[ T^1\Aut^0_{W,d}=B^0(\Hom^*(W,W))=\{df+fd\mid f\in\Hom^{-1}(W,W)\}.\]
Since all the functors are smooth and homogeneous, it is sufficient
to prove that
the exponential induces  isomorphisms on the tangent spaces.
\end{proof}

The considerations made in Section~\ref{sec.coarsegrass} lead us to 
define  the functors
\[
M_{V,W},\Grass_{V,W}\colon\mathbf{Art}\to \mathbf{Sets}\]
\[ M_{V,W}(A)=\{f\in \Aut_W(A)\mid df(V\otimes A)\subseteq f(V\otimes
A)\},\]
\[ \Grass_{V,W}=\frac{M_{V,W}}{\Aut^0_{W,d}\times \Aut_{V,W}}.\]

\begin{proposition}\label{prop.defundgrass}
In the notation above, let $\chi\colon L_{V,W}\to L_W$ the inclusion.
Then there exist natural isomorphisms of functors
\[ \MC_{\chi}\mapor{\sim}M_{V,W},\qquad \Def_{\chi}\mapor{\sim}\Grass_{V,W}.\]
\end{proposition}

\begin{proof}

Note that, since $V$ is a subcomplex of $(W,d)$, we have $d(V\otimes
A)\subseteq V\otimes A$, and so
\[
M_{V,W}(A)=\{f\in \Aut_W(A)\mid (f^{-1}df-d)(V\otimes A)\subseteq
V\otimes A\}=\]
\[=\{f\in \Aut_W(A)\mid f^{-1}df-d\in
L_{V,W}^1\otimes\mathfrak{m}_A\},\]
and then, by Lemma~\ref{lem.expo} and by the identity
$e^{-a}\ast 0=e^{-a}de^a-d$,
\[ M_{V,W}(A)=\{e^a\in \exp(L^0_W\otimes\mathfrak{m}_A)\mid
e^{-a}\ast 0\in
L_{V,W}^1\otimes\mathfrak{m}_A\}.\]
Recall that
\[\MC_{\chi}(A)
=\!\left\{(x,e^a)\in (L^1_{V,W}\otimes\mathfrak{m}_A)\times
\exp(L^0_W\otimes\mathfrak{m}_A)\mid
dx+\frac{1}{2}[x,x]=0,\;e^{-a}\ast0=\chi(x)\right\}.\]
Since $\chi$ is injective and
the set of solutions of the Maurer-Cartan equation in
$L^1_W\otimes\mathfrak{m}_A$ is preserved by the gauge action
$l\mapsto e^{-a}\ast l$, we have
\[\MC_{\chi}(A)
=\left\{e^a\in
\exp(L^0_W\otimes\mathfrak{m}_A)\mid
e^{-a}\ast0\in L^1_{V,W}\otimes\mathfrak{m}_A\right\}\]
and then the isomorphism $\MC_{\chi}\simeq M_{V,W}$.\\
The gauge action
\[(\exp(L_{V,W}^0\otimes\mathfrak{m}_A)\times
\exp(B^0(L_W\otimes\mathfrak{m}_A)))
\times \MC_{\chi}(A)\mapor{\ast}\MC_{\chi}(A)\]%
becomes
\[ (e^l, e^{dm})\ast e^a= e^{dm}e^ae^{-\chi(l)}
=e^{dm}e^ae^{-l}.\]
According to Lemma~\ref{lem.expo} we have
\[\exp(L_{V,W}^0\otimes\mathfrak{m}_A)\times
\exp(B^0(L_W\otimes\mathfrak{m}_A))=\Aut_{V,W}(A)\times \Aut^0_{W,d}(A)\]
and it is immediate to observe that the gauge action
is identified with the natural group action
\[
\Aut_{V,W}(A)\times \Aut^0_{W,d}(A)\times M_{V,W}(A)\to  M_{V,W}(A)
\]
and then we have the isomorphism of quotient functors
\[
\Grass_{V,W}=\Def_\chi.
\]
\end{proof}

Note that the isomorphism $\Def_\chi\to \Grass_{V,W}$
has a natural explicit description: it is the map
\[e^a \mapsto (e^a(V\otimes A),d)\subseteq (W\otimes A,d).\]
\par

\begin{theorem}\label{thm.pregrassinjectivecohomology}
Let $V$ be a subcomplex of $(W,d)$ such that the inclusion $V\hookrightarrow W$
induces an injective morphism in cohomology
$H^*(V)\hookrightarrow H^*(W)$.
Identifying $H^*(V)$ with its image in
$H^*(W)$, then the cohomology functor $H^*$ gives a natural
transformation of functors
\[ H^*\colon \Grass_{V,W}\to
\Grass_{H^*(V),H^*(W)}\simeq
\prod_{i}\Grass_{H^i(V),H^i(W)}.\]
\end{theorem}

\begin{proof}
Recall that
\[ \Grass_{H^*(V),H^*(W)}(A)=
\{\text{free $A$ submodules }F_A\subseteq H^*(W)\otimes A\mid
F_A\otimes_A\K=H^*(V)\}.\]
Given $e^a\in \MC_{\chi}(A)$ denote by $d_a=e^{-a}de^a\colon
  W\otimes
A\to W\otimes A$; then $(V\otimes A,d_a)$ is a subcomplex of
$(W\otimes A,d_a)$ and
$e^a$ is a morphism of complexes
$(W\otimes A,d_a)\to (W\otimes A,d)$. By local flatness
criteria, the cohomology of $(V\otimes A,d_a)$ is a free
$A$-module and the map
\[
H^*(e^a(V\otimes A),d)\to
H^*(W\otimes A,d)\simeq H^*(W)\otimes A
\]
is injective since it factors as
\[ H^*(e^a(V\otimes A),d)\simeq
  H^*(V\otimes A,d_a)\xrightarrow{H^*(e^a)}
H^*(W\otimes A,d)\]
where the isomorphism on the left is $H^*(e^{-a})$.
Therefore the map
\[ {h}\colon \MC_{\chi}(A)\to \Grass_{H^*(V),H^*(W)}(A),\]
\[{h}(e^a)=\text{image of the natural map }H^*(e^a(V\otimes A))
\to H^*(W)\otimes A\]
is well defined and factors to a natural transformation of functors
\[ H^*\colon \Grass_{V,W}\mapor{}
\Grass_{H^*(V),H^*(W)}.\]
\end{proof}

\bigskip

\section{Homotopy invariance}

\begin{lemma}\label{lem.qitoqi} If $V\hookrightarrow W$
is a
quasiisomorphism, then the inclusion
\[ L_{V,W}\hookrightarrow L_W.\]
and the projection
\[ L_{V,W}\to L_V,\qquad f\mapsto f\vert_V\]
are quasiisomorphisms of DGLA.
\end{lemma}

\begin{proof}
We have two short exact sequences
\[
0\to L_{V,W}\to L_W \to \Hom^*(V,W/V)\to 0
\]
\[
0\to \Hom^*(W/V,W)\to L_{V,W}\to L_V \to 0
\]
Since the complex $W/V$ is acyclic and
  the
bifunctor $\Hom^*$ commutes with cohomology,
the complexes $\Hom^*(V,W/V)$ and $\Hom^*(V,W/V)$ are acyclic.
\end{proof}

\begin{lemma}\label{lem.croce}
Assume we have an exact diagram of differential graded vector
spaces
\[\begin{array}{ccccccccc}
&&&&0&&&&\\
&&&&\mapver{}&&&&\\
&&&&A&&&&\\
&&&&\mapver{\alpha}&&&&\\
0&\mapor{}&B&\mapor{\varphi}&C&\mapor{\psi}&D&\mapor{}&0\\
&&&&\mapver{\beta}&&&&\\
&&&&E&&&&\\
&&&&\mapver{}&&&&\\
&&&&0&&&&\end{array}\]%
If $\psi\alpha$ is a quasiisomorphism, then also $\beta\varphi$ is a
quasiisomorphism.
\end{lemma}

\begin{proof}
For every $i$ in $\mathbb{Z}$ the map $H^i(\psi\alpha)\colon H^i(A)\to
H^i(D)$ is an isomorphism. Since $H^i(\psi\alpha)=H^i(\psi)\circ
H^i(\alpha)$, the map $H^i(\alpha)\colon H^i(A)\to H^i(C)$
is injective and $H^i(\psi)\colon H^i(C)\to H^i(D)$ is surjective.
Therefore, the long exact sequences
\[
\cdots\to H^{i-1}(E)\xrightarrow{\delta^{i-1}}
H^{i}(A)\xrightarrow{H^{i}(\alpha)} H^{i}(C)
\xrightarrow{H^{i}(\beta)}
H^{i}(E)\xrightarrow{\delta^i} H^{i+1}(A)\to \cdots
\]
and
\[
\cdots\to H^{i-1}(D)\xrightarrow{\delta^{i-1}}
H^{i}(B)\xrightarrow{H^{i}(\varphi)} H^{i}(C)
\xrightarrow{H^{i}(\psi)}
H^{i}(D)\xrightarrow{\delta^i} H^{i+1}(B)\to \cdots
\]
can be refined as
\[
H^{i-1}(E)\xrightarrow{\delta^{i-1}}0\to
H^{i}(A)\xrightarrow{H^{i}(\alpha)}
H^{i}(C)\xrightarrow{H^{i}(\beta)}
H^{i}(E)\xrightarrow{\delta^i}0
\]
and
\[
H^{i-1}(D)\xrightarrow{\delta^{i-1}}0\to
H^{i}(B)\xrightarrow{H^{i}(\varphi)} H^{i}(C)\xrightarrow{H^{i}(\psi)}
H^{i}(D)\xrightarrow{\delta^i}0
\]
and
we have an
exact diagram
\[\begin{array}{ccccccccc}
&&&&0&&&&\\
&&&&\mapver{}&&&&\\
&&&&H^i(A)&&&&\\
&&&&\mapver{H^i(\alpha)}&&&&\\
0&\mapor{}&H^i(B)&\mapor{H^i(\varphi)}&H^i(C)&\mapor{H^i(\psi)}&H^i(D)
&\mapor{}&0\\
&&&&\mapver{H^i(\beta)}&&&&\\
&&&&H^i(E)&&&&\\
&&&&\mapver{}&&&&\\
&&&&0&&&&\end{array}\]%

Now it is an easy exercise on linear algebra to show that
$H^i(\beta\varphi)=H^i(\beta)\circ H^i(\varphi)$ is an isomorphism,
i.e.,
$\beta\varphi$ is a quasiisomorphism. Note that, by symmetry, if
$\beta\varphi$ is a quasiisomorphism, then also
$\psi\alpha$ is a quasiisomorphism.
\end{proof}

\begin{lemma}\label{lem.terna}
Consider a subcomplex $U$ of $V\subseteq W$
and denote
\[ L_{U,V,W}=L_{U,W}\cap L_{V,W}=\{f\in\Hom^*(W,W)\mid
f(U)\subseteq U,\;  f(V)\subseteq V\}.\]%
\begin{enumerate}

\item If $U\hookrightarrow V$ is a quasiisomorphism, then the two inclusions
\[ L_{U,V,W}\hookrightarrow L_{V,W},\qquad L_{U,V,W}\hookrightarrow
L_{U,W}\] are quasiisomorphisms of DGLA.

\item If $V\hookrightarrow W$ is a quasiisomorphism, then the inclusions
$L_{U,V,W}
\hookrightarrow
L_{U,W}$ and the projection
\[ L_{U,V,W}\to L_{U,V},\qquad f\mapsto f_{|V}\]
are quasiisomorphisms of DGLA.
\end{enumerate}

\end{lemma}

\begin{proof}
Assume that $U\hookrightarrow V$ is a quasiisomorphism;  consider the 
exact sequences
\[ 0\to L_{U,V,W}\xrightarrow{({\rm id}\vert_V,{\rm id})}
L_{U,V}\oplus
L_{V,W}
\xrightarrow{\tau}L_V\to 0,\qquad
\tau(f,g)=f-g_{|V}\]
and
\[ 0\to L_{U,V}\xrightarrow{({\rm id},0)} L_{U,V}\oplus L_{V,W}
\xrightarrow{\pi_2} L_{V,W}\to 0.\]
The composition $\tau\circ({\rm
id},0)\colon L_{U,V}\to L_V$ coincides with the inclusion
$L_{U,V}\hookrightarrow
L_V$.
By Lemma
\ref{lem.qitoqi}, the inclusion
$L_{U,V}\hookrightarrow
L_V$ is a quasiisomorphism.  Therefore, by
Lemma
\ref{lem.croce} the map $\pi_2\circ({\rm id}\vert_V,{\rm
id})\colon\allowbreak L_{U,V,W}\to L_{V,W}$ is a quasiisomorphism,
i.e.,
$L_{U,V,W}\hookrightarrow L_{V,W}$ is a quasiisomorphism.
\par
The acyclicity of $V/U$, together with the exact sequence

\[  0\to L_{U,V,W}\to L_{U,W}
\xrightarrow{\sigma}\Hom^*\left(\frac{V}{U},\frac{W}{V}\right)\to
0,\qquad
\sigma(g)=g_{|V},\]
imply that $L_{U,V,W}\hookrightarrow L_{U,W}$
is a quasiisomorphism.\\
The proof of the second statement is very similar.
Assume that $V\hookrightarrow W$ is a quasiisomorphism and
consider the exact sequences
\[ 0\to L_{U,V,W}\xrightarrow{({\rm id}\vert_V,{\rm id})}
L_{U,V}\oplus
L_{V,W}
\xrightarrow{\tau}L_V\to 0,\qquad
\tau(f,g)=g_{|V}-f\]
and
\[ 0\to L_{V,W}\xrightarrow{(0,{\rm id})} L_{U,V}\oplus L_{V,W}
\xrightarrow{\pi_1} L_{U,V}\to 0.\]
The projection
$L_{V,W}\to
L_V$ is a quasiisomorphism by Lemma \ref{lem.qitoqi}.
This projection can be written as the composition $\tau\circ(0,{\rm
id})\colon L_{V,W}\to L_V$, hence by
Lemma
\ref{lem.croce} the map $\pi_1\circ({\rm id}\vert_V,{\rm
id})\colon\allowbreak L_{U,V,W}\to L_{U,V}$ is a quasiisomorphism,
i.e.,
the projection
$L_{U,V,W}\to L_{U,V}$ is a quasiisomorphism. Moreover $W/V$ is
acyclic and then the inclusion $L_{U,V,W}\hookrightarrow L_{U,W}$ is a
quasiisomorphism.
\end{proof}

\begin{lemma}\label{lem.homotopyinvarianceofgrass}
Assume we have a commutative diagram of \textbf{inclusions}
of differential graded
vector spaces
\[ \begin{array}{ccc}
U&\rightarrow&V\\
\downarrow&&\downarrow\\
Z&\rightarrow&W\end{array}\]%
If the horizontal arrows are
quasiisomorphisms then $\Grass_{U,Z}\simeq \Grass_{U,W}\simeq
\Grass_{V,W}$.
\end{lemma}

\begin{proof}
In the notation of Lemma \ref{lem.terna},
there exists a commutative diagram of DGLA, where every arrow is
either the natural inclusion or the natural
projection:
\[\begin{array}{ccccccccc}
L_{U,Z}&\!\leftarrow\!&L_{U,Z,W}&\!\rightarrow\!&L_{U,W}
&\!\leftarrow\!&L_{U,V,W}&\!\rightarrow\!&L_{V,W}\\
\mapver{\eta}&&\mapver{}&&\mapver{\rho}&&\mapver{}&&\mapver{\chi}\\
L_Z&\!\leftarrow\!&L_{Z,W}&\!\rightarrow\!&L_W
&\!\leftarrow\!&L_W&\!\rightarrow\!&L_W
\end{array}\]
By  Lemma \ref{lem.terna} the horizontal arrows are quasiisomorphism and then,
Theorem~\ref{thm.basic} gives an isomorphism of deformation functors
\[
\Grass_{U,Z}\simeq \Def_{\eta}\simeq
\Def_{\rho}\simeq\Grass_{U,W}
\simeq \Def_{\chi}\simeq\Grass_{V,W}.
\]
\end{proof}

\begin{remark}\label{rem.cartesian}
In the same hypothesis of Lemma~\ref{lem.homotopyinvarianceofgrass},
if in addition $U=V\cap Z$, then the isomorphism
$\Grass_{U,Z}\simeq \Grass_{V,W}$ can be done in a more explicit and easy way.
In fact the two inclusion
\[ \alpha \colon L_{Z,W}\cap L_{V,W}\hookrightarrow L_{U,V,W},\qquad
\beta \colon L_{Z,W}\cap L_{V,W}\hookrightarrow L_{U,Z,W}\]
are quasiisomorphisms since their cokernels are the acyclic complexes
\[
\operatorname{Coker}(\alpha)=\Hom^*\left(\frac{Z}{U},\frac{W}{Z}\right),\qquad
\operatorname{Coker}(\beta)=\Hom^*\left(\frac{V}{U},\frac{W}{V}\right).\]
Then, according to Lemma~\ref{lem.homotopyinvarianceofgrass}
every horizontal arrow of
\[\begin{array}{ccccc}
L_{U,Z}&\!\leftarrow\!&L_{Z,W}\cap
L_{V,W}&\!\hookrightarrow\! &L_{V,W}\\
\mapver{\eta}&&\mapver{}&&\mapver{\chi}\\
L_Z&\!\leftarrow\!&L_{Z,W}&\!\hookrightarrow\!&L_W
\end{array}\]
is a quasiisomorphism. The isomorphism
$\Grass_{V,W}\simeq \Grass_{U,Z}$ can then be explicitly
described as follows: given $A\in\mathbf{Art}$ and an element 
$e^a(V\otimes A)$ in
$\Grass_{V,W}(A)$, there exists (and it is unique up to gauge)
an element $e^\alpha\in
\exp(\Hom^0(W,W)\otimes\mathfrak{m}_A)$ such that
$e^\alpha(V\otimes A)=e^a(V\otimes A)$ and $e^\alpha(Z\otimes 
A)\subseteq Z\otimes A$. Then the
isomorphism $\Grass_{V,W}\simeq \Grass_{U,Z}$ is given by
the map
\[
e^a(V\otimes A) \mapsto e^\alpha(U\otimes A)\subseteq Z\otimes A.
\]
\end{remark}

\begin{theorem}\label{thm.grassinjectivecohomology}
Let $V$ be a subcomplex of $(W,d)$ such that the inclusion $V\hookrightarrow W$
induces an injective morphism in cohomology
$H^*(V)\hookrightarrow H^*(W)$.
Identifying $H^*(V)$ with its image in
$H^*(W)$, then the cohomology functor $H^*$ gives an isomorphism of functors
\[ H^*\colon \Grass_{V,W}\mapor{\sim}
\Grass_{H^*(V),H^*(W)}\cong \prod_{i}\Grass_{H^i(V),H^i(W)}.\]
\end{theorem}

\begin{proof}
It is possible to find ``harmonic representatives''
${\mathcal H}_V\subseteq V$ and ${\mathcal H}_W\subseteq W$
with
${\mathcal H}_V={\mathcal H}_W\cap V$. Indeed, the
injectivity of
$H^*(V)\to H^*(W)$ implies the two equalities
\[ Z^*(V)=V\cap Z^*(W),\qquad B^*(V)=V\cap B^*(W).\]
Therefore it is possible to find a  subspace
$\mathcal{H}_W\subseteq Z^*(W)$ such that
$\mathcal{H}_W\oplus B^*(W)=Z^*(W)$ and
$(\mathcal{H}_W\cap V)\oplus B^*(V)=Z^*(V)$. Set
${\mathcal H}_V=\mathcal{H}_W\cap V$. Then the diagram
\[ \begin{array}{ccc}
\mathcal{H}_V&\mapor{}&V\\
\mapver{}&&\mapver{}\\
\mathcal{H}_W&\mapor{}&W\end{array}\]%
satisfies the hypothesis of
Lemma~\ref{lem.homotopyinvarianceofgrass} and Remark~\ref{rem.cartesian}
and we have an isomorphism
$\Grass_{V,W}\xrightarrow{\sim}\Grass_{\mathcal{H}_V,\mathcal{H}_W}$ given
by
\[
e^a(V\otimes A) \mapsto e^\alpha({\mathcal H}_V\otimes A),
\]
where $e^\alpha\in \exp(\Hom^0(W,W)\otimes {\mathfrak m}_A)$
is such that $e^\alpha(V\otimes A)=e^a(V\otimes A)$ and $e^\alpha({\mathcal
H}_W\otimes A)\subseteq {\mathcal H}_W\otimes A$.\\
Now, according to
Theorem~\ref{thm.pregrassinjectivecohomology}
it is sufficient to note that the natural
map
\[e^\alpha({\mathcal H}_V\otimes A)\to H^*(e^\alpha(V\otimes A),d)=
H^*(e^a(V\otimes A),d)\]
is an isomorphism.
\end{proof}

\bigskip
\part{The universal period map}

In this part we work over the field $\K=\mathbb{C}$ of complex numbers.
Unless otherwise specified, the symbol $\otimes$ denotes tensor
product over $\mathbb{C}$.

\bigskip
\section{The Kodaira-Spencer DGLA}
\label{sec.kodairaspencer}

We will follows the same general notation of \cite{Voisin}; in particular, for
a differentiable manifold $M$  we denote by
$\mathcal{A}_M^{p}$ the sheaf of differentiable $p$-forms with
complex coefficients
and by $d\colon \mathcal{A}_M^{p}\to \mathcal{A}_M^{p+1}$ the De Rham
differential.\\
We think an almost complex structure on $M$ as a  subsheaf
$\mathcal{V}\subseteq \mathcal{A}_M^{1}$ of locally free
$\mathcal{A}_M^{0}$-modules such that
$\mathcal{V}\oplus \overline{\mathcal{V}}=\mathcal{A}_M^{1}$.\\

For a complex manifold $X$ we denote by:\begin{itemize}

\item $T_{X,\mathbb{C}}=T^{1,0}_X\oplus T^{0,1}_X$ the complexified
differential tangent
bundle.

\item $T_X\simeq T^{1,0}_X$ the holomorphic tangent bundle.

\item $\mathcal{A}_X^{p,q}$ the sheaf of differentiable $(p,q)$-forms
and by $\mathcal{A}_X^{p,q}(T_X)$ the sheaf of $(p,q)$-forms
with values in $T_X$.

\item ${A}_X^{p,q}$  and ${A}_X^{p,q}(T_X)$ the vector spaces of
global sections
of $\mathcal{A}_X^{p,q}$ and $\mathcal{A}_X^{p,q}(T_X)$ respectively.

\end{itemize}

Recall that an almost complex structure $\mathcal{V}\subseteq
\mathcal{A}_X^{1}$
is called integrable if there exists a
structure of complex manifold on $X$ such that
$\mathcal{V}=\mathcal{A}_X^{1,0}$.\\

The direct sum
\[\mathcal{A}_X=\mathop{\oplus}_i\mathcal{A}^{i}_X,\quad\text{ where }\quad
\mathcal{A}^{i}_X=\mathop{\oplus}_{p+q=i}\mathcal{A}^{p,q}_X,\]
endowed with the wedge product $\wedge$, is a sheaf of graded
algebras; we denote by
$\DER^{a,b}(\mathcal{A}_X)$ the sheaf of its $\mathbb{C}$-linear
derivations of bidegree $(a,b)$.
Notice that $\de$ and $\debar$ are global sections of
$\DER^{1,0}(\mathcal{A}_X)$ and $\DER^{0,1}(\mathcal{A}_X)$ respectively.\\
The direct sum
$\DER^{*}(\mathcal{A}_X)=\oplus_k\oplus_{a+b=k}\DER^{a,b}(\mathcal{A}_X)$
is a sheaf of differential graded Lie algebras, with
its natural bracket
\[ [f,g]=fg-(-1)^{\deg(f)\deg(g)}gf\]
and  differential $[d,-]=[\de+\debar,-]$.\\
Similarly $\mathcal{A}^{0,*}_X(T_X)$ is a sheaf of DGLA, where  the
bracket $[\;,\;]$ and the differential
$D$ are defined in local coordinates by the formulas:
\[ D\left(\phi\frac{\de~}{\de z_i}\right)=-\debar(\phi)\frac{\de~}{\de
z_i},\qquad \phi\in \mathcal{A}_X^{0,*},\]
\[\left[fd\bar{z}_{I}\desude{~}{z_{i}},
gd\bar{z}_{J}\desude{~}{z_{j}}\right]=
d\bar{z}_{I}\wedge
d\bar{z}_{J}\left(f\desude{g}{z_{i}}\desude{~}{z_{j}}-
g\desude{f}{z_{j}}\desude{~}{z_{i}}\right).\]
The contraction of differential forms with vector fields is used to
define two
injective morphisms of sheaves:
the \emph{contraction map}
\[ \bi\colon \mathcal{A}^{0,*}_X(T_X)\to
\DER^{*}(\mathcal{A}_X)[-1],\qquad
a\mapsto \bi_{a},\quad \bi_{a}(\omega)=a\contr \omega,\]%
and the \emph{holomorphic Lie derivative}
\[ \bl\colon \mathcal{A}^{0,*}_X(T_X)\to
\DER^{*}(\mathcal{A}_X),\qquad
a\mapsto \bl_{a}=[\de,\bi_a],\quad \bl_{a}(\omega)=\de(a\contr
\omega)+(-1)^{\deg(a)}a\contr \de\omega.\]%

\begin{lemma}\label{lem.cartanequalities}
In the notation above,
for every $a,b\in \mathcal{A}^{0,*}_X(T_X)$  we have
\[\bi_{Da}=-[\debar, \bi_a],\qquad
\bi_{[a,b]}=[\bi_a,[\de,\bi_b]],\qquad
[\bi_a,\bi_b]=0.\]
In particular, since $[\de,\bi_b]=[d,\bi_b]+\bi_{Db}$, the contraction
map $\bi$ is a Cartan homotopy and
the holomorphic Lie derivative is a morphism of sheaves of differential graded
Lie algebras.
\end{lemma}

\begin{proof}
Since locally $\mathcal{A}_X$ is generated as
$\mathbb{C}$-algebra by
$\mathcal{A}_X^{0,0}\oplus\mathcal{A}_X^{1,0}\oplus\mathcal{A}_X^{0,1}$,
if
$\min(a,b)\le -2$ then $\DER^{a,b}(\mathcal{A}_X)=0$ , while if
$\min(a,b)\le -1$ then every $h\in \DER^{a,b}(\mathcal{A}_X)$ is
$C^{\infty}$-linear.
The elements of the third formula belong to
$\DER^{-2,*}(\mathcal{A}_X)$
and therefore vanish. Every term of first two formulas belongs to
$\DER^{-1,*}(\mathcal{A}_X)$ and then it is sufficient to check
equalities of such derivations when applied to $dz_i$, where
$z_1,\ldots,z_n$ are local holomorphic coordinates.
This is straightforward and it is left to
the reader: see also Lemma 7 of \cite{CCK}.\end{proof}

The contraction map gives a natural isomorphism of vector spaces
\[\bi\colon A^{0,1}_X(T_X)\to
\Hom_{\mathcal{A}^{0}_X}(\mathcal{A}^{1,0}_X,\mathcal{A}^{0,1}_X).\]

For every sufficiently small $\xi\in A^{0,1}_X(T_X)$, the graph of
$\bi_{\xi}\in
\Hom_{\mathcal{A}^{0}_X}(\mathcal{A}^{1,0}_X,\mathcal{A}^{0,1}_X)$
determines a variation of the almost complex structure of $X$ given by
\[ \mathcal{A}^{1,0}_{\xi}=\{\omega\in\mathcal{A}^{1}_X\mid
\pibar(\omega)=\bi_{\xi}\pi(\omega)\}=\{\omega\in\mathcal{A}^{1}_X\mid
\pibar(\omega)=\xi\contr\pi(\omega)\},\]%
where $\pi\colon \mathcal{A}^{1}_X\to \mathcal{A}^{1,0}_X$ and
$\pibar\colon \mathcal{A}^{1}_X\to \mathcal{A}^{0,1}_X$ are the projections.
Then denote
\[\mathcal{A}^{0,1}_{\xi}=\overline{\mathcal{A}^{1,0}_{\xi}},
\qquad \mathcal{A}^{p,q}_{\xi}=
\bigwedge^p\mathcal{A}^{1,0}_{\xi}\otimes \bigwedge^q\mathcal{A}^{0,1}_{\xi}.\]
The sheaf of $\xi$-holomorphic functions is by definition
\[ \mathcal{O}_{\xi}=\{f\in \mathcal{A}^{0}_X\mid df\in
\mathcal{A}^{1,0}_{\xi}\}\]
and then, according to the definition of $\mathcal{A}^{1,0}_{\xi}$, we have
\[ \mathcal{O}_{\xi}=\{f\in \mathcal{A}^{0}_X\mid \debar f=\xi\contr\de f
\}=\{f\in \mathcal{A}^0_X\mid (\debar+\bl_{\xi})f=0\}.\]

Since $\bi_{\xi}$ is a nilpotent derivation of degree 0 of
the sheaf of graded algebras $\mathcal{A}_X$, its exponential
\[ e^{\bi_{\xi}}=\exp(\xi\contr)\colon \mathcal{A}_X\to \mathcal{A}_X\]
is an isomorphism of graded algebras.

\begin{lemma} In the above notation $e^{\bi_{\xi}}(\mathcal{A}^{1,0}_X)=
\mathcal{A}^{1,0}_{\xi}$.
\end{lemma}

\begin{proof}
For every $\omega\in \mathcal{A}^1_X$ we have
$e^{\bi_{\xi}}(\omega)=\omega+\xi\contr\pi(\omega)$, therefore
\[\pi e^{\bi_{\xi}}(\omega)=\pi(\omega),\qquad
\pibar e^{\bi_{\xi}}(\omega)=\pibar(\omega)+\xi\contr\pi(\omega)\]
and then $e^{\bi_{\xi}}(\omega)\in \mathcal{A}^{1,0}_{\xi}$
if and only if $\pibar(\omega)=0$.
\end{proof}

The Newlander-Nirenberg theorem \cite{NN}, \cite[Thm. 5.5]{Kobook}
implies (see e.g. \cite[Lecture 1]{Catacime}) that the following
four conditions are equivalent:
\begin{enumerate}

\item The almost complex structure $\mathcal{A}^{1,0}_{\xi}$ is integrable.

\item $\xi$ is a solution of the Maurer-Cartan equation
\cite[Equation 5.86]{Kobook}, \cite[Equation 2.5]{Catacime}
\[ D\xi+\frac{1}{2}[\xi,\xi]=0.\]

\item For every $x\in X$ there exist
$f_1,\ldots,f_n\in \mathcal{O}_{\xi,x}$, $n=\dim X$, such that
$df_1,\ldots, df_n$ are a basis of the $\mathcal{A}^0_{X,x}$-module
$\mathcal{A}^{1,0}_{\xi,x}$.

\item  $d\mathcal{F}^1_{\xi}\subseteq \mathcal{F}^1_{\xi}$, where
$\mathcal{F}^1_{\xi}=\oplus_{p\ge 1}\mathcal{A}^{p,q}_{\xi}=
(\mathcal{A}^{1,0}_{\xi})\subseteq \mathcal{A}_{X}$
is the graded ideal sheaf generated by  $\mathcal{A}^{1,0}_{\xi}$.
\end{enumerate}

\begin{definition}
The Kodaira-Spencer algebra of a complex manifold $X$
is the differential graded Lie
algebra  $K_X=\oplus_{i}A^{0,i}_X(T_X)$ of global sections of
$\oplus_{i}\mathcal{A}^{0,i}_X(T_X)$.
\end{definition}

Denoting by $\Def_X\colon\mathbf{Art}\to\mathbf{Set}$
the functor of infinitesimal deformations of $X$:
\[ \Def_X(B)=\frac{\text{deformations of $X$ over $\Spec(B)$}}{\sim}.\]
A deformation of $X$ over $\Spec(B)$ may be interpreted as
a morphism
$\mathcal{O}_B\to \mathcal{O}_X$ of sheaves of
$B$-algebras such that $\mathcal{O}_B$ is flat over $B$ and
the induced map $\mathcal{O}_B\otimes_B\mathbb{C}\to \mathcal{O}_X$
is an isomorphism.\\

The following result is well known \cite{clemens}, \cite{GoMil2},
\cite[Ex. 3.4.1]{K}:
a detailed proof will also appear in
\cite{Iaconophd}

\begin{theorem}
There exists an isomorphism of functors
\[ \mathcal{O}\colon \Def_{K_X}\to \Def_X\]
defined in the following way: given a local Artinian
$\mathbb{C}$-algebra $B$ and a solution of the
Maurer-Cartan equation
$\xi\in A^{0,1}_X(T_X)\otimes\mathfrak{m}_B$ we set
\[\mathcal{O}_{\xi}=\ker(\mathcal{A}^{0,0}_X\otimes
B\xrightarrow{\debar+\bl_{\xi}}
\mathcal{A}^{0,1}_X\otimes B)=\{f\in \mathcal{A}^{0,0}_X\otimes B\mid
\debar f=\xi\contr\de f\}\]
and the map  $\mathcal{O}_{\xi}\to \mathcal{O}_{X}$ is induced by the
projection
$\mathcal{A}^{0,0}_X\otimes B\to \mathcal{A}^{0,0}_X\otimes
\mathbb{C}=\mathcal{A}^{0,0}_X$.
\end{theorem}

\bigskip
\section{The period map}

Let $X$ be a fixed complex manifold: we shall denote by
$A_X=F^0\supseteq F^1\supseteq\cdots$ the Hodge filtration of
differential forms on $X$, i.e. for every $p\ge 0$
\[F^p=\mathop{\bigoplus}_{i\ge p}\mathop{\bigoplus}_{j}A^{i,j}_X.\]

\begin{theorem}\label{thm.periodasLinfinity}
Let $p$ be a fixed nonnegative integer and
consider the inclusion of differential graded Lie algebras
\[  L_{F^p,A_X}=\{f\in \Hom^*(A_X,A_X)\mid f(F^p)\subseteq F^p\}
\mapor{\chi} \Hom^*(A_X,A_X)=L_{A_X}.\]
Then the linear map
\[ \mathfrak{p}^p\colon K_X\to C_{\chi},\qquad
\mathfrak{p}^p(\xi)=(\bl_{\xi},\bi_{\xi})=([\de,\bi_{\xi}],\bi_{\xi})\]%
is a linear $L_{\infty}$-morphism.\\
In particular $\mathfrak{p}^p$ induces a natural transformation of functors:
\[ \mathfrak{p}^p\colon \Def_X\to \Def_{\chi}=\Grass_{F^p,A_X}.\]
\end{theorem}

\begin{proof}
According to Lemma~\ref{lem.cartanequalities} and
Proposition~\ref{prop.cartan},
the Lie derivative
\[ K_X=A^{0,*}_X(T_X)\mapor{\bl}\Hom^*(A_X,A_X)\]
is a morphism of differential graded Lie algebras
and the map
\[ \tilde{\bi}\colon K_X\to C_{\bl}=K_X\oplus\Hom^*(A_X,A_X)[-1],\qquad
\tilde{\bi}(\xi)=(\xi,\bi_{\xi})\]%
is a linear $L_{\infty}$-morphism.\\
The morphism $\mathfrak{p}^p$ is the composition of $\tilde{\bi}$ and
the linear $L_{\infty}$-morphism
$C_{\bl}\to C_{\chi}$ induced
by the horizontal arrows of the following
commutative diagram of differential graded Lie algebras
\[ \begin{array}{ccc}
A^{0,*}_X(T_X)&\mapor{\bl}&L_{F^p,A_X}\\
\mapver{\bl}&&\mapver{\chi}\\
\Hom^*(A_X,A_X)&\mapor{=}&L_{A_X}\end{array}\]
\end{proof}

\begin{corollary}\label{cor.newdifferential}
If $\xi\in A^{0,1}_X(T_X)$ is a solution of the Maurer-Cartan
equation, then in the associative
algebra $\Hom^*(A_X,A_X)$ we have the equality
\[ e^{-\bi_{\xi}} d e^{\bi_{\xi}}=
d+e^{-\bi_{\xi}}\ast 0=d+\bl_{\xi}=d+[\de,\bi_{\xi}].\]
\end{corollary}

\begin{proof}
The first equality follows from the explicit description of the gauge
action in the
DGLA $\Hom^*(A_X,A_X)$. Taking $p=0$ in
Theorem~\ref{thm.periodasLinfinity}, the pair
$\mathfrak{p}^0(\xi)=(\bl_{\xi},\bi_{\xi})$
satisfies the Maurer-Cartan equation in
$C_{\chi}$ and then $e^{\bi_{\xi}}\ast \bl_{\xi}=0$.\\
Obviously the above equality can be also proved directly as a consequence of
Cartan homotopy formulas (Lemma~\ref{lem.cartanequalities}).
\end{proof}

\begin{theorem}\label{thm.periodsforkaehler}
In the same notation of Theorem~\ref{thm.periodasLinfinity},
assume that $X$ is a compact K\"{a}hler manifold and denote
by
\[ \mathcal{P}^p\colon \Def_X\to \Grass_{H^*(F^p), H^*(A_X)}\]
the composition of the natural transformation
$\mathfrak{p}^p\colon \Def_X\to \Grass_{F^p,A_X}$
and the cohomology isomorphism
$H^*\colon \Grass_{F^p,A_X}\to\Grass_{H^*(F^p),H^*(A_X)}$ (see
Theorem~\ref{thm.grassinjectivecohomology}). Then $\mathcal{P}^p$ is
the universal period map.
\end{theorem}

\begin{proof}
Let $B\in\mathbf{Art}$ and $\xi\in\MC_{K_X}(B)$; by definition
\[ \mathcal{P}^p(\xi)=H^*(e^{\bi_{\xi}}(F^p\otimes B))\subseteq
H^*(A_X\otimes B)=H^*(A_X)\otimes B.\]
On the other hand, the period of the infinitesimal deformation
$\mathcal{O}_{\xi}=\ker(\debar+\bl_{\xi})$ is the $B$-submodule
$H^*(F^p_{\xi})\subseteq H^*(A_X\otimes B)$, where $F^p_{\xi}$ is the
complex of global sections
of the differential ideal sheaf $\mathcal{F}^p_{\xi}\subseteq
\mathcal{A}_X\otimes B$ generated
by
$(d\mathcal{O}_{\xi})^p$.\\
It is sufficient to prove that $e^{\bi_{\xi}}(F^p\otimes B)=F^p_{\xi}$;
since $e^{\bi_{\xi}}\colon \mathcal{A}_X\otimes B\to 
\mathcal{A}_X\otimes B$ is a
morphism of sheaves of $B$-algebras,
it is sufficient to prove that $e^{-\bi_{\xi}}(d\mathcal{O}_{\xi})\subseteq
\mathcal{A}_X^{1,0}\otimes B$. This equality, together rank
considerations, will imply that
\[e^{-\bi_{\xi}}(\mathcal{F}^p_{\xi})=\mathop{\bigoplus}_{i\ge
p,j}\mathcal{A}_X^{i,j}\otimes B.\]
Since $e^{\bi_{\xi}}$ is the identity on $\mathcal{A}_X^{0,0}\otimes B$, by
Corollary~\ref{cor.newdifferential} we can write
\[
e^{-\bi_{\xi}}(d\mathcal{O}_{\xi})=e^{-\bi_{\xi}}de^{\bi_{\xi}}\mathcal{O}_{\xi}=
(\de+\debar+\bl_{\xi})\mathcal{O}_{\xi}=
\de \mathcal{O}_{\xi}\subseteq \de\mathcal{A}_X^{0,0}\otimes B
\subseteq \mathcal{A}_X^{1,0}\otimes B.\]
\end{proof}

\begin{remark}
In the statement of Theorem~\ref{thm.periodsforkaehler} the K\"{a}hler
assumption is used in cohomological sense; more precisely we only
require that the cohomology of the complex $F^p$ injects into the De
Rham cohomology of $X$.
\end{remark}

\begin{corollary}[Griffiths]\label{cor.griffiths}
The differential of the universal period map is
\[d\mathcal{P}^p=\bi\colon H^1(X,T_X)\to
\bigoplus_i \Hom\left(F^pH^i(X,\mathbb{C}),
\frac{H^{i}(X,\mathbb{C})}{F^pH^{i}(X,\mathbb{C})}\right).\]
\end{corollary}

\begin{proof}

Recall from Remark~\ref{remark.conequotient} that the projection on the second
factor induces an identification
\[
H^j(C_\chi)\xrightarrow{\sim}
H^{j-1}(\Hom^*(F^p,A_X/F^p))=\bigoplus_i\Hom(H^i(F^p),H^{i+j-1}(A_X/F^p)).
\]
Via this identification, the $L_\infty$-morphism $\mathfrak{p}^p\colon
K_X\to C_\chi$, induces in cohomology the map
\[
H^j(\mathfrak{p}^p)=\bi\colon H^j(K_X)\to
\bigoplus_i\Hom(H^i(F^p),H^{i+j-1}(A_X/F^p)).
\]
The differential of ${\mathcal P}^p\colon \Def_X\to
\Grass_{F^p,A_X}\simeq \Grass_{H^*(F^p),H^*(X;{\mathbb C})}$ is
therefore
\[
d{\mathcal P}^p=\bi\colon H^1(K_X)\to
\bigoplus_i\Hom(H^i(F^p),H^{i}(A_X/F^p)).
\]

\end{proof}

\begin{corollary}[Kodaira's Principle \cite{clemens,CCK,ranUVHS}]
The obstructions to
deformations of $X$ are contained in the kernel of
\[ \bi\colon H^2(X,T_X)\to
\bigoplus_i \Hom\left(F^pH^i(X,\mathbb{C}),
\frac{H^{i+1}(X,\mathbb{C})}{F^pH^{i+1}(X,\mathbb{C})}\right),\]
for every $p\ge 0$.
\end{corollary}
\begin{proof}
We use the same general argument of \cite[Section 5]{CCK}: since the period map
$\mathcal{P}^p\colon \Def_X \to \Grass_{H^*(F^p),H^*(X,\mathbb{C})}$ is induced
by the $L_\infty$-morphism ${\mathfrak p}^p\colon K_X\to C_\chi$, the
linear map
\[
H^2({\mathcal P}^p)\colon H^2(X,T_X)\to H^2(C_\chi)
\]
is a morphism of obstruction theories, i.e., it commutes with the
natural obstruction maps for $\Def_X$ and
$\Grass_{H^*(F^p),H^*(X,\mathbb{C})}$
\cite{semireg}.  In particular, obstructions to deforming the complex
structure of the K\"ahler manifold $X$ are mapped to obstructions of
the functor $\Grass_{H^*(F^p),H^*(X,\mathbb{C})}$.  Since the latter
is unobstructed,  the obstructions to deforming
$X$ are annihilated by $H^2(\mathfrak{p}^p)$.  By the proof of
Corollary~\ref{cor.griffiths}, $H^2(\mathfrak{p}^p)=\bi$.
\end{proof}

\bigskip
\section{Trasversality}

Consider a fixed compact K\"{a}hler manifold $X$ and a differential 
graded commutative
$\mathbb{C}$-algebra $(\Omega, d_{\Omega})$. Let
\[  \operatorname{Id}\otimes d_{\Omega}\colon A_X\otimes \Omega\to 
A_X\otimes \Omega,\qquad
(\operatorname{Id}\otimes d_{\Omega})(a\otimes
\omega)= (-1)^{\deg(a)}a\otimes d_{\Omega}(\omega)\]
the trivial extension of  $d_{\Omega}$.
Then $\operatorname{Id}\otimes d_{\Omega}$ is a differential of the graded
algebra $A_X\otimes \Omega$ inducing
a flat connection
\[  \operatorname{Id}\otimes d_{\Omega}\colon
H^*(A_X)\otimes \Omega^i\to H^*(A_X)\otimes \Omega^{i+1}.\]
Assume now that $B$ is a commutative unital $\mathbb{C}$-algebra and let
$\phi\colon B\to \Omega^0$ be a morphism of graded unital
$\mathbb{C}$-algebras.\\
Via the natural isomorphisms
\[A_X\otimes \Omega=(A_X\otimes B)\otimes_B \Omega,\qquad
H^*(A_X)\otimes \Omega=H^*(A_X\otimes B)\otimes_B \Omega\]
the operator $\operatorname{Id}\otimes d_{\Omega}$ induces the
differential
\[  \nabla\colon (A_X\otimes B)\otimes_B \Omega\to (A_X\otimes 
B)\otimes_B \Omega,\]
\[\nabla((a\otimes b)\otimes_B \omega)=
(-1)^{\deg(a)}(a\otimes 1)\otimes_B d_{\Omega}(\phi(b)\omega)
\]
and the flat connection
\[  \nabla\colon H^*(A_X\otimes B)\otimes_B \Omega^i\to 
H^*(A_X\otimes B)\otimes_B \Omega^{i+1},\]
that, by analogy with the case of $B$ a power series ring and
$\Omega=\bigwedge_B^*\Omega^1_B$ we shall call \emph{Gauss-Manin
connection}.\\

Assume now that $B\in\mathbf{Art}$ and consider
a deformation of $X$ over
$\operatorname{Spec}(B)$
determined by the Kuranishi data $\xi\in \MC_{K_X}(B)\subseteq
A^{0,1}_X(T_X)\otimes\mathfrak{m}_B$.

The classical Griffiths' trasversality theorem \cite{Greencime},
\cite[Prop. 10.18]{Voisin} generalizes to the following result:

\begin{proposition}[Trasversality]\label{prop.trasversality}
Let $F^p_{\xi}$ be the
complex of global sections
of the differential ideal sheaf $\mathcal{F}^p_{\xi}\subseteq
\mathcal{A}_X\otimes B$ generated
by $(d\mathcal{O}_{\xi})^p$. Then
\[ \nabla (H^*(F^p_{\xi})\otimes_B \Omega^i)\subseteq
H^*(F^{p-1}_{\xi})\otimes_B \Omega^{i+1}.\]
\end{proposition}

\begin{proof}
It is not restrictive to assume $B\subseteq \Omega^0\subseteq \Omega$
and $\phi$ the inclusion.  For notational simplicity denote by
$\delta\colon A_X^{0,*}(T_X)\otimes \Omega\to A_X^{0,*}(T_X)\otimes
\Omega$ the trivial extension of $d_\Omega$, i.e. $\delta(\xi\otimes
b)=(-1)^{\deg(\xi)}\xi\otimes d_\Omega(b)$.\\
It is sufficient to prove that $\nabla (F^p_{\xi}\otimes_B
\Omega^i)\subseteq F^{p-1}_{\xi}\otimes_B \Omega^{i+1}$.\\

The contraction map $\bi$ extends naturally to a Cartan homotopy
(Example~\ref{ex.cambiobasepercartan}) \[ \bi\colon
A_X^{0,*}(T_X)\otimes \Omega\to \Hom^*(A_X,A_X)\otimes \Omega\subseteq
\Hom^*((A_X\otimes B)\otimes_B \Omega,(A_X\otimes B)\otimes_B\Omega)\]
such that $\bi_{\delta x}=-[\nabla,\bi_x]$ ; in particular
$[\bi_{x},[\bi_{x},\nabla]]=-[\bi_{x},\bi_{\delta x}]=0$ for every
$x$.  For every nilpotent $a\in \Hom^0((A_X\otimes B)\otimes_B
\Omega,(A_X\otimes B)\otimes_B\Omega)$ we have the equality
\begin{multline*}
\qquad\qquad [\nabla,e^a]=\nabla e^a-e^a\nabla=e^a(e^{-a}\nabla
e^a-\nabla)\\
=e^a(e^{-\ad_{a}}(\nabla)-\nabla) =e^a[\nabla,a]-e^a\sum_{n\ge
2}\frac{(-\ad_{a})^n}{n!}(\nabla).\qquad\qquad
\end{multline*}
Since $\xi\in \MC_{K_X}(B)$, for $a\in A_X\otimes B$ and $\omega\in
\Omega$ we have
$\bi_{\xi}(a\otimes_B\omega)=\bi_{\xi}(a)\otimes_B\omega$; therefore
\[ [\nabla,e^{\bi_{\xi}}]=e^{\bi_{\xi}}\left(-\bi_{\delta\xi}\pm
\frac{1}{2}[\bi_{\xi},\bi_{\delta\xi}]
\pm\frac{1}{6}[\bi_{\xi},[\bi_{\xi},\bi_{\delta\xi}]]\pm\cdots\right)=
-e^{\bi_{\xi}}\bi_{\delta\xi}.\] Now it is easy to conclude the proof:
we have seen in the proof of Theorem~\ref{thm.periodsforkaehler} that
$F^p_{\xi}=e^{\bi_{\xi}}(F^p\otimes B)$ and then every $v\in
F^p_{\xi}\otimes_B \Omega$ can be written as $v=e^{\bi_{\xi}}(u)$,
with $u\in (F^p\otimes B)\otimes_B \Omega$.  Therefore
\begin{multline*}
\nabla(v)=\nabla(e^{\bi_{\xi}}(u))=[\nabla,e^{\bi_{\xi}}](u)+e^{\bi_{\xi}}(\nabla
u)\\
=e^{\bi_{\xi}}(-\bi_{\delta\xi}(u)+\nabla(u))\in
e^{\bi_{\xi}}((F^{p-1}\otimes B)\otimes_B \Omega)=
F^{p-1}_{\xi}\otimes_B \Omega.
\end{multline*}
\end{proof}

\end{document}